\documentclass[11pt, amssymb]{amsart}
\usepackage{amsmath, amsthm, amsfonts, verbatim, amssymb}
\usepackage{eucal}

\DeclareFontFamily{OT1}{rsfs}{}
\usepackage[all]{xy}
\CompileMatrices

\DeclareFontShape{OT1}{rsfs}{n}{it}{<-> rsfs10}{}
\DeclareMathAlphabet{\mathscr}{OT1}{rsfs}{n}{it}
\def\sG{{\mathscr G}}

\def\fF{\mathfrak F}

\def\sM{{\mathscr M}}
\def\oF{\overline F}

\def\osP{\overline{\mathscr P}}

\def\osG{\overline{\mathscr{G}}}
\def\osT{\overline{\mathscr{T}}}
\def\orG{\overline{G}^{\,\rm pred}}
\def\cB{\mathscr B}

\def\sR{{\mathscr R}}

\def\sV{{\mathscr V}}

\def\bC{{\mathbb C}}
\def\bR{{\mathbb R}}
\def\bZ{{\mathbb Z}}

\def\({{\Big(}}
\def\){{\Big)}}

\def\cB{{\mathcal B}}

\def\sS{{\mathscr S}}
\def\cG{{\mathcal G}}
\def\cS{{\mathcal S}}
\def\sT{{\mathscr T}}

\def\sZ{{\mathscr Z}}

\def\cC{{\mathcal C}}

\def\ok{\overline k}

\def\X{\widetilde X}

\def\ni{\noindent}

\def\ok{\overline k}

\def\oS{\overline S}

\def\ok{\overline k}

\def\osH{\overline{\mathscr{H}}}
\def\oC{\overline C}

\def\oF{\overline F}
\def\oC{\overline C}
\def\sV{\mathcal V}

\def\wK{\widehat{K}}
\def\wk{\widehat{k}}
\def\fK{\mathfrak{K}}
\def\osS{{\overline{\mathscr{S}}}}
\def\osT{{\overline{\mathscr{T}}}}
\def\X{{\mbox{X}}}
\def\fo{\mathfrak{o}}
\def\fp{\mathfrak{p}}
\def\fT{\mathfrak{T}}

\def\fZ{\mathfrak{Z}}

\def\fG{\mathfrak{G}}

\def\cO{\mathcal{O}}
\def\cA{\mathcal{A}}
\def\cF{\mathcal{F}}
\def\cC{\mathcal{C}}
\def\sU{\mathscr{U}}
\def\osU{\overline{\mathscr{U}}}
\def\osZ{\overline{\mathscr{Z}}}
\def\sV{\mathscr{V}}

\def\sB{\mathscr{B}}
\def\sD{\mathscr{D}}
\def\sN{\mathscr{N}}
\def\cG{\mathcal{G}}
\def\cP{\mathcal{P}}
\def\sI{\mathscr{I}}

\def\cS{\mathcal{S}}
\def\cX{\mathcal{X}}
\def\cY{\mathcal{Y}}

\begin{document}
\title[Unramified descent in Bruhat-Tits theory]
{A new approach to unramified descent in Bruhat-Tits theory}
 \maketitle
\vskip-7mm
{\centerline {{By} {\sc {Gopal Prasad}}}

\vskip4mm
\centerline{\it Dedicated to my wife Indu Prasad with gratitude}
\vskip5mm

\begin{abstract} We present a new approach to unramified (or ``\'etale'') descent in Bruhat-Tits theory of reductive groups over a discretely valued field $k$ with Henselian valuation ring, which appears to be conceptually simpler, and more geometric, than the original approach of  Bruhat and Tits. We are able to derive the main results of the theory over $k$ from the theory over the maximal unramified extension $K$ of $k$. Even in the most interesting  case for number theory and representation theory, where $k$ is a locally compact nonarchimedean field, the geometric approach described in this paper appears to be considerably simpler than the original approach. 
\end{abstract} 
\vskip4mm

\ni{\bf Introduction.}  Let $k$ be a field given with a nontrivial $\bR$-valued nonarchimedean valuation $\omega$. We will assume throughout this paper that the valuation  ring   $$\fo := \{x\in k^{\times}\, \vert\, \omega(x)\geqslant 0\}\cup \{0\}$$ of this valuation is  Henselian. Let $K$ be the maximal unramified extension of $k$ (the term ``unramified extension'' includes the condition that the residue field extension is separable, so the residue field of $K$ is the separable closure $\kappa_s$ of the residue field $\kappa$ of $k$).  While discussing Bruhat-Tits theory in this introduction, and beginning with 1.6 
 everywhere in the paper, we will assume that $\omega$ is a discrete valuation. Bruhat-Tits theory of connected reductive $k$-groups $G$ that are quasi-split over $K$ (i.e., $G_K$ contains a Borel subgroup, or, equivalently, the centralizer of a maximal $K$-split torus of $G_K$ is a torus) has two parts. The first part of the theory, developed in [BrT2,\,\S4], is called the ``quasi-split descent'' (``descente quasi-d\'eploy\'ee'' in French), and is due to Iwahori-Matsumoto, Hijikata, and Bruhat-Tits. It is the theory over $K$ assuming that $G_K$ is quasi-split.  If $G_K$ is quasi-split, the group $G(K)$ has a rather simple structure. In particular, it admits a ``Chevalley-Steinberg system'', which is used in [BrT2,\,\S4] to get a valuation on root datum, that in turn is used to construct the Bruhat-Tits building $\cB(G/K)$ of $G(K)$. In [BrT2, Intro.] Bruhat and Tits say that ``La descente quasi-d\'eploy\'ee est la plus facile.''. The second part, called the ``unramified descent'' (``descente \'etale'' in French),  is due to Bruhat and Tits. This part derives  Bruhat-Tits theory of $G$ (over $k$), and also the Bruhat-Tits building of $G(k)$, from Bruhat-Tits theory of $G$ over $K$ and the Bruhat-Tits building of $G(K)$, using descent of valuation of root datum from $K$ to $k$. This second part is quite technical; see [BrT1,\,\S9] and [BrT2,\,\S5].  
\vskip1mm

The purpose of this paper is to present an alternative approach to unramified descent which appears to be conceptually simpler, and more geometric,  than the approach in [BrT1], [BrT2], in that it does not use descent of valuation of root datum from $K$ to $k$ to show that $\cB(G/K)^{\Gamma}$, where $\Gamma$ is the Galois group of $K/k$, is an affine building.  In this approach, we will use Bruhat-Tits theory, and the buildings, only over the maximal unramified extension $K$ of $k$ and derive the main results of the theory for reductive groups over $k$. In \S4, we discuss the notions of hyperspecial points and hyperspecial parahoric subgroups and describe conditions under which these exist. In \S5, we define a natural filtration of the root groups $U_a(k)$ and also introduce a valuation of the root datum of $G/k$ relative to a maximal $k$-split torus $S$, using the geometric results of \S\S2,\,3 that provide the Bruhat-Tits building of $G(k)$. The approach described here appears to be considerably simpler than the original approach even for reductive groups over locally compact nonarchimedean fields\:(i.e., discretely valued complete fields with finite residue field). In \S6, we prove results over discretely valued fields with Henselian valuation ring and perfect residue field of dimension $\leqslant 1$; of these, Theorems 6.1 and 6.2 may be new.     

\vskip3mm

\ni{\it Acknowledgements.} Lemma 1.3 is well-known; I thank Brian Conrad for writing its  
proof given below and for his comments on several earlier versions of this paper. I thank 
Vladimir Chernousov, Philippe Gille, Tasho Kaletha, Dipendra Prasad, Igor Rapinchuk, Mark Reeder, 
Bertrand R\'emy, Guy Rousseau and Jiu-Kang Yu for their comments, corrections  and suggestions. I thank the referee for very carefully 
reading the paper and for her/his helpful comments. 
I was partially supported by NSF-grant DMS-1401380.
\vskip3mm

\ni{\bf 1. Preliminaries} 
\vskip3mm

We assume below that $k$ is a field given with a nontrivial $\bR$-valued nonarchimedean valuation $\omega$ and the valuation  ring $\fo$ of $\omega$ is Henselian. It is known that the valuation ring $\fo$ is Henselian if and only if the valuation $\omega$ extends uniquely to any algebraic extension of $k$. The valuation ring of a discretely valued {complete} field is always Henselian. (For various equivalent definitions of the Henselian property, see [Be,\,\S\S2.3-2.4].)   Let $K$ be the maximal unramified extension of $k$ contained in a fixed algebraic closure $\ok$ of $k$. We will denote the unique valuation on the algebraic closure $\ok\, (\supset K)$, extending the given valuation on $k$, also by $\omega$.  The residue field of $k$ will be denoted by $\kappa$ and the valuation ring of $K$ by $\cO$. The residue field of $K$ is the separable closure $\kappa_s$ of $\kappa$. We will denote the Galois group  ${\rm{Gal}}(K/k) = {\rm{Gal}}({\kappa_s}/\kappa)$ by $\Gamma$. For a set $Z$ given with an action of $\Gamma$, we will denote by $Z^{\Gamma}$ the subset consisting of elements fixed under $\Gamma$. 

\vskip1mm

\ni{\it Bounded subsets.} Let $X$ be an affine $k$-variety.  A subset $B$ of $X(\ok)$ is said to be {\it bounded} if for every $f\in k[X]$, the set $\{\omega(f(b))\:\vert\: b\in B\}$ is bounded below.  If $X$ is an affine algebraic $k$-group and $B$, $B'$ are two nonempty subsets of $X(\ok)$, then $BB' = \{bb'\:\vert\: b\in B, b'\in B'\}$ is bounded if and only if both $B$ and $B'$ are bounded. 
\vskip1mm

Let $G$ be a connected reductive group defined over $k$.  The group of $k$-rational characters on $G$ will be denoted by ${\rm{X}}_k^{\ast}(G)$.   The following theorem is due to Bruhat, Tits and Rousseau. An elementary proof was given in [P] 
which we recall here for the reader's convenience. 
\vskip2mm

\ni{\bf Theorem 1.1.} {\em $G(k)$ is bounded if and only if $G$ is anisotropic over $k$}.
\vskip2mm

\ni{\bf Remark 1.2.} Thus if $k$ is a nondiscrete locally compact field, then $G(k)$ is compact if and only if $G$ is $k$-anisotropic. 

\vskip1mm

We fix a faithful $k$-rational representation of $G$ on a finite dimensional $k$-vector space $V$ and view $G$ as a $k$-subgroup of ${\rm GL}(V)$. To prove the above theorem we will use the following two lemmas: 

\vskip1mm

\ni{\bf Lemma 1.3.} {\em If $f:X\rightarrow Y$ is a finite $\ok$-morphism between affine $\ok$-schemes of finite type and $B$ is a bounded subset of $Y(\ok)$, then the subset $f^{-1}(B)$ of $X(\ok)$ is bounded.}
\vskip1mm

\ni{\it Proof.} Since $\ok[X]$ is module-finite over $\ok[Y]$, we can pick a finite set of generators of $\ok[X]$ as a $\ok[Y]$-module (so also as a $\ok[Y]$-algebra), and each satisfies a monic polynomial over $\ok[Y]$.   Hence, this realizes $X$ as a closed 
subscheme of the closed subscheme $Z \subset Y \times {\mathbb A}^n$ defined by $n$ monic 1-variable polynomials $f_1(t_1), \ldots\,,  f_n(t_n)$ over $\ok[Y]$, so it remains to observe that when one has a bound on the coefficients of a monic 1-variable polynomial over $\ok$ of known degree (e.g., specializing any $f_j$ at a $\ok$-point of $Y$) then one gets a bound on its possible $\ok$-rational roots depending only on the given coefficient bound and the degree of the monic polynomial.   \hfill$\Box$
\vskip2mm

\ni{\bf Lemma 1.4.}  {\em Let $\cG$ be an unbounded subgroup of $G(k)$ which is dense in $G$ in the Zariski-topology. Then $\cG$ contains an element $g$ which has an eigenvalue $\alpha$ (for the action on $V$) with $\omega(\alpha)<0$.}
\vskip1mm

\ni{\it Proof.} Let $$\ok\otimes_k V =: V_0\supset V_1\supset\cdots\, \supset V_s\supset V_{s+1}=\{0\}$$ be a flag of $G_{\ok}$-invariant subspaces such that for $0\leqslant i\leqslant s$, the natural representation $\varrho_i$ of $G_{\ok}$ on $W_i :=V_i/V_{i+1}$ is irreducible. Let $\varrho =\bigoplus_i \varrho_i$ be the representation of $G_{\ok}$ on $\bigoplus_i W_i$. The kernel of $\varrho$ is obviously a unipotent normal $\ok$-subgroup scheme of the reductive group $G_{\ok}$, and hence it is finite. Now as $\cG$ is an unbounded subgroup of $G(k)$, Lemma 1.3 implies that $\varrho(\cG)$ is an unbounded subgroup of $\varrho(G(\ok))$.  Hence, there is a non-negative integer $a\leqslant s$ such that $\varrho_a(\cG)$ is unbounded. 
\vskip1mm

 Since $W_a$ is an irreducible $G_{\ok}$-module, and $\cG$ is dense in $G$ in the Zariski-toplogy, $\varrho_a(\cG)$ spans ${\rm End}_{\ok}(W_a)$. We fix $\{g_i\}\subset \cG$ so that 
$\{\varrho_a(g_i)\}$ is a basis of  ${\rm End}_{\ok}(W_a)$. Let $\{ f_i\}\subset {\rm End}_{\ok}(W_a)$ be the basis which is dual to the basis $\{\varrho_a(g_i)\}$ with respect to the trace-form. Then ${\rm Tr}(f_i\cdot\varrho_a(g_j))=\delta_{ij}$, where $\delta_{ij}$ is the Kronecker's delta. Now assume that the eigenvalues of all the elements of $\cG$ lie in the valuation ring $\fo_{\ok}$ of $\ok$. Then for all $x\in \cG$, ${\rm Tr}(\varrho_a(x))$ is contained in $\fo_{\ok}$.  For $g\in \cG$, if $\varrho_a(g) = \sum c_if_i$, with $c_i\in \ok$, then ${\rm Tr}(\varrho_a(g\cdot g_j)) = \sum_i c_i{\rm Tr}(f_i\cdot \varrho_a(g_j)) = c_j$. As  ${\rm Tr}(\varrho_a(g\cdot g_j))\in \fo_{\ok}$, we conclude that $c_j$ belongs to the ring of integers $\fo_{\ok}$ for all $j$ (and all $g\in \cG$). This implies that $\varrho_a(\cG)$ is bounded, a contradiction. \hfill$\Box$
\vskip2mm

\ni{\it Proof of Theorem 1.1.} As ${\rm GL}_1(k) =k^{\times}$ is unbounded, we see that if $G$ is $k$-isotropic, then $G(k)$ is unbounded. We will now assume that $G(k)$ is unbounded and prove the converse.  

\vskip1mm

It is well known that $G(k)$ is dense in $G$ in the Zariski-topology [Bo,\:18.3], hence according to Lemma 1.4, there is an element $g\in G(k)$ which has an eigenvalue $\alpha$ with $\omega(\alpha)\ne 0$. Now in case $k$ is of positive characteristic, after replacing $g$ by a suitable positive integral power, we assume that $g$ is semi-simple. On the other hand, in case $k$ is of characteristic zero, let $g =s\cdot u=u\cdot s$ be the Jordan decomposition of $g$ with $s\in G(k)$ semi-simple and $u\in G(k)$ unipotent. Then the eigenvalues of $g$ are same as that of $s$. So, after replacing $g$ with $s$, we may (and do) again assume that $g$ is semi-simple.  There is a maximal $k$-torus $T$ of $G$ such that $g\in T(k)$ (see [BoT], Proposition 10.3 and Theorem 2.14(a); note that according to Theorem 11.10 of [Bo], $g$ is contained in a maximal torus of $G$).  Since any absolutely irreducible representation of a torus is $1$-dimensional, there exists a finite Galois extension $\fK$ of $k$ and a character $\chi$ of $T_{\fK}$ such that $\chi(g) = \alpha$.  Then $$\omega\big(\big(\sum_{\gamma\in {\rm Gal}(\fK/k)} {^\gamma{\chi}}\big)(g)\big) = m\omega(\chi(g))= m\omega(\alpha)\ne 0;$$ where $m = [\fK:k]$. Thus the character $ \sum_{\gamma\in {\rm Gal}(\fK/k)} {^\gamma{\chi}}$ is nontrivial. On the other hand, this character is obviously defined over $k$. Hence, $T$ admits a nontrivial character defined over $k$ and therefore it contains a nontrivial $k$-split subtorus. This proves that if $G(k)$ is unbounded, then $G$ is isotropic over $k$. \hfill$\Box$  
\vskip2mm

\ni{\bf Proposition 1.5.} {\em We assume that the derived subgroup $G':=(G,G)$ of $G$ is $k$-anisotropic. Then $G(k)$ contains a unique maximal bounded subgroup $G(k)_b$; it has the following description: $$G(k)_b = \{\, g\in G(k)\ \vert\  \chi(g)\in \fo^{\times} \ {\rm {for \ all}}\  {\chi \in}\, {\rm{X}}_k^{\ast}(G)\}. $$}
\ni{\it Proof.}  Let $G_a$ be the inverse image of the maximal $k$-anisotropic subtorus of the $k$-torus $G/G'$ under the natural homomorphism $G\rightarrow G/G'$. Then $G_a$ be the maximal connected normal $k$-anisotropic subgroup of $G$.  Let $S$ be the maximal $k$-split central torus of $G$. Then $G = S\cdot G_a$ (almost direct product). Let $C = S\cap G_a$; $C$ is a finite central $k$-subgroup scheme, so $G_a/C$ is $k$-anisotropic. Let $f: G\rightarrow G/C =(S/C)\times (G_a/C)$ be the natural homomorphism. The image of the induced homomorphism $f^{\ast}\,:\, \X_k^{\ast}((S/C)\times (G_a/C))\rightarrow \X_k^{\ast}(G)$  is of finite index. It is obvious that as $(G_a/C)(k)$ is bounded (by Theorem 1.1), the proposition is true for the direct product $(S/C)\times (G_a/C)$. Now using Lemma 1.3 we conclude that the proposition holds for $G$. \hfill$\Box$    

\vskip2mm

\ni{\em We shall henceforth assume that the valuation $\omega$ on $k$ is discrete.}

\vskip2mm

\ni {\bf 1.6.} Let $S$ be a maximal $k$-split torus of $G$, $Z(S)$ its centralizer in $G$ and $Z(S)' = (Z(S),Z(S))$ the derived subgroup of $Z(S)$. Then $Z(S)'$ is a connected semi-simple group which is anisotropic over $k$ since $S$ is a maximal $k$-split torus of $G$. Hence, by Theorem 1.1,  $Z(S)'(k)$ is bounded, and according to  Proposition 1.5,  $Z(S)(k)$ contains a unique maximal bounded subgroup $Z(S)(k)_b$. This maximal bounded subgroup admits the following description:  $$Z(S)(k)_b = \{ z\in Z(S)(k)\ | \ \chi(z)\in \fo^{\times} \ \ {\rm for\ all }\  \ \chi\in {\X_k^{\ast}}(Z(S))\}.$$
\vskip1mm

The restriction map ${\X_k^{\ast}}(Z(S))\rightarrow {\X_k^{\ast}}(S)$ is injective and its image is of finite index in ${\X_k^{\ast}}(S)$. Let ${\X_{\ast}}(S)= {\rm{Hom}}_k ({\rm{GL}}_1, S)$ and $V(S) =\bR\otimes_{\bZ} {\X_{\ast}}(S)$.  Let the homomorphism $\nu: Z(S)(k)\rightarrow V(S)$ be defined by: 
$$\chi(\nu(z)) = -\omega(\chi(z)) \ \  {\rm{for}}  \ z\in Z(S)(k) \ {\rm{and}} \  \chi\in {\X_k^{\ast}}(Z(S))\,(\hookrightarrow {\X_k^{\ast}}(S)).$$ 
Then $Z(S)(k)_b$ is the kernel of $\nu$. As the image of $\nu$ is isomorphic to $\bZ^r$, $r =\dim S$, we conclude that $Z(S)(k)/Z(S)(k)_b$ is isomorphic to $\bZ^r$.
\vskip2mm 
 
 \ni{\bf 1.7. Fields of dimension $\leqslant 1$ and a theorem of Steinberg.} A field $\fF$ is said to be of dimension $\leqslant 1$ if 
 finite dimensional central simple algebras with center a finite separable extension of $\fF$ are matrix algebras [S, Ch.\,II, \S3.1]. For example, every finite field is of dimension $\leqslant 1$.  
 \vskip1mm
 
 We now recall the following theorem of Steinberg: {\em For a smooth connected linear algebraic group $\fG$ defined over a field $\fF$ of dimension $\leqslant 1$, the Galois cohomology ${\rm{H}}^1(\fF, \fG)$ is trivial if either $\fF$ is perfect or $\fG$ is reductive}\:[S, Ch.\,III, Thm.\,$1'$ and Remark (1) in \S2.3]. This vanishing theorem implies that {\em if $\fF$ is of dimension $\leqslant 1$, then a connected linear algebraic $\fF$-group $\fG$ is quasi-split, i.e., it contains a Borel subgroup defined over $\fF$}, assuming that $\fG$ is reductive when $\fF$ is not perfect; note that the proof of the fact that assertion $({\rm{i}}')$ of Theorem 1 in [S, Ch.\,III,\,\S2.2] implies that the semi-simple group $L$ contains a Borel subgroup defined over the base field does not require the base field to be perfect.    
\vskip1mm

We assume in this paragraph that the residue field $\kappa$ of $k$ is perfect. Then the residue field of the maximal unramified extension $K$ is the algebraic closure $\overline\kappa$ of $\kappa$.  Let  $\wK$ denote the completion of $K$. The discrete valuation on $K$ extends uniquely to the completion $\wK$ and the residue field $\overline\kappa$ of $K$ is also the residue field of $\wK$. Hence, by Lang's theorem, $\wK$ is a $({\rm{C}}_1)$-field [S, Ch.\,II, Example 3.3(c) in \S3.3], so it is of dimension $\leqslant 1$\,[S, Ch.\,II, Corollary in \S3.2]. According to a well-known result  (see, for example, Proposition 3.5.3(2) of [GGM] whose proof simplifies considerably in the smooth affine case), for any smooth algebraic $K$-group $\fG$, the natural map ${\rm{H}}^1(K,\fG)\rightarrow {\rm{H}}^1(\wK, \fG)$ is bijective. This result, combined with the above theorem of Steinberg, implies that  for any connected reductive $K$-group $\fG$, ${\rm{H}}^1(K,\fG)$ is trivial, hence every connected reductive $K$-group is quasi-split.  
\vskip2mm

\ni{\bf Notation.} Given a smooth connected linear algebraic group $\fG$ defined over a field $\fF$, 
we will denote its $\fF$-unipotent radical, i.e., the maximal smooth connected normal unipotent $\fF$-subgroup, 
by $\sR_{u,\fF}(\fG)$. The quotient $\fG^{\rm{pred}}:=\fG/\sR_{u,\fF}(\fG)$ is {\it pseudo-reductive}{\footnote{For definition and properties of pseudo-reductive groups and pseudo-parabolic subgroups, see [CGP] or [CP].}}; it is the maximal pseudo-reductive quotient of $\fG$. If the field $\fF$ is perfect, then pseudo-reductive groups are actually reductive and pseudo-parabolic subgroups of $\fG$ are parabolic subgroups. 
  
\vskip2mm
  
 Let $\fZ$ be the maximal $k$-torus contained in the center of $G$ that splits over $K$. There is a natural action of the  Galois group $\Gamma$ of $K/k$ on ${\rm{Hom}}_K({\rm{GL}}_1, \fZ_K)$ and ${\rm{Hom}}_K({\rm{GL}}_1, \fZ_K)^{\Gamma} = {\rm{Hom}}_k({\rm{GL}}_1, \fZ)$. Let $V(\fZ_K)=\bR\otimes_{\bZ}{\rm{Hom}}_K({\rm{GL}}_1, \fZ_K)$. The  action of $\Gamma$ on ${\rm{Hom}}_K({\rm{GL}}_1, \fZ_K)$ extends to an $\bR$-linear action on  $V(\fZ_K)$, and $V(\fZ_K)^{\Gamma} = \bR\otimes_{\bZ}{\rm{Hom}}_k({\rm{GL}}_1, \fZ)$.  We will denote the derived subgroup $(G,G)$ of $G$ by $G'$  throughout this paper. $G'$ is the maximal connected normal semi-simple subgroup of $G$ and there is a natural bijective correspondence between the set of maximal $K$-split tori of $G'_K$ and the set of maximal $K$-split tori of $G_K$ given by $T\mapsto \fZ_KT$. 
 \vskip2mm

 {\em In the rest of the paper we will assume that  Bruhat-Tits theory is available for $G$ over $K$,  that is, there is an affine building $\cB(G/K)$, called the ``enlarged'' Bruhat-Tits building of $G(K)$, on which this group acts by isometries, and given a nonempty bounded subset $\Omega$ of an apartment of this building, there is a smooth affine $\cO$-group scheme  $\sG^{\circ}_{\Omega}$ with generic fiber $G$ and connected special fiber--the building $\cB(G/K)$ and the group schemes $\sG^{\circ}_{\Omega}$ having the properties described in 1.8 and 1.9 below. }
 \vskip2mm
 
   When $G$ is quasi-split over $K$--for example, if the residue field $\kappa$ of $k$ is perfect\,(1.7)--then Bruhat-Tits theory is available for $G$ over $K$; see [BrT2,\,\S4]. 
\vskip2mm
 
  \ni {\bf 1.8. Bruhat-Tits theory for $G$ over $K$:}  There exists an affine building called the Bruhat-Tits building of $G(K)$. It carries a $G(K)$-invariant metric and a natural structure of a polysimplicial complex on which $G(K)$ acts by polysimplicial automorphisms. This building is also the Bruhat-Tits building of $G'(K)$, and we will denote it by $\cB(G'/K)$. As we are assuming that $\omega$ is a discrete valuation, $\cB(G'/K)$ is complete\,[BrT1, Thm.\,2.5.12(i)].  The apartments of this building are in bijective correspondence with maximal $K$-split tori of $G_K$.  If $A$ is the apartment of $\cB(G'/K)$ corresponding to a maximal $K$-split torus $T$ of $G_K$, then for $g\in G(K)$, $g\cdot A$ is also an apartment and it corresponds to the maximal $K$-split torus $gTg^{-1}$. Hence the stabilizer of $A$ in $G(K)$ is $N(T)(K)$, where $N(T)$ denotes the normalizer of $T$ in $G_K$. 
\vskip.5mm

The facets of $\cB(G'/K)$ of maximal dimension are called {\it chambers}.  The group $G'(K)$ acts transitively on the set of ordered pairs consisting of an apartment of $\cB(G'/K)$ and a chamber in it. In particular, $N(T)(K)\cap G'(K)$ acts transitively on the set of chambers in $A$. 
\vskip.5mm 

There is a natural action of $G(K)$ on the Euclidean space $V(\fZ_K)$ by translations, with $G'(K)$ acting trivially.  The direct product $V(\fZ_K)\times \cB(G'/K)$ carries a $G(K)$-invariant metric extending the metric on $\cB(G'/K)$, and there is visibly an action of $V(\fZ_K)$ on this product by translation in the first factor. This direct product is called the {\it enlarged} Bruhat-Tits building of $G(K)$ and we will denote it by $\cB(G/K)$. The apartments of the enlarged building are by definition the subspaces of the form $V(\fZ_K)\times A$, where $A$ is an apartment of $\cB(G'/K)$.     If $G$ is semi-simple, i.e., $G' = G$, then $\cB(G/K) = \cB(G'/K)$.
\vskip.5mm

Let $T$ be a maximal $K$-split torus of $G_K$ and $A$ be the corresponding apartment of $\cB(G/K)$. Then $A$ is an affine space under  $V(T):=\bR\otimes_{\bZ}{{\rm X}_{\ast}}(T)$, where ${\rm X}_{\ast}(T)= {\rm{Hom}}_K({\rm{GL}}_1, T)$,  and $N(T)(K)$ acts on $A$ by affine transformations which we will describe now. Let ${\rm Aff}(A)$ be the group of affine automorphisms of $A$ and $\nu: N(T)(K)\rightarrow {\rm Aff}(A)$ be the action map. For $n\in N(T)(K)$, the derivative $d\nu(n): V(T)\rightarrow V(T)$ is the map induced by the action of $n$ on ${\rm{X}}_{\ast}(T)$ (i.e.,\,the Weyl group action).    So for $z\in Z(T)(K)$, $d\nu(z)$ is the identity, hence $\nu(z)$ is a translation; this translation is described by the following formula: 
$$\chi(\nu(z)) = -\omega(\chi(z)) \ {\rm for \ all}\  \chi \in {{\rm X}_K^{\ast}}(Z(T))\,(\hookrightarrow {{\rm X}_K^{\ast}}(T)),$$
here we regard the translation $\nu(z)$ as an element of  $V(T)$. Since for $z$ in the maximal bounded subgroup $Z(T)(K)_b$ of $Z(T)(K)$, $\omega(\chi(z)) = 0$ for all $\chi\in  {{\rm X}_K^{\ast}}(Z(T))$ (Proposition 1.5), the above formula shows  that  $Z(T)(K)_b$ acts trivially on $A$.  
\vskip.5mm

Given two points $x$ and $y$ of $\cB(G/K)$, there is a unique geodesic $[xy]$ joining them and this geodesic lies in every apartment which contains $x$ and $y$. A subset of the building is called {\it convex} if for any $x,y$ in the set, the geodesic $[xy]$ is contained in the set. For a subset $X$ of $\cB(G/K)$, $\overline{X}$ will denote its closure. If $X$ is convex, then so is $\overline{X}$.

\vskip1mm

 Let $G(K)^{\natural}$ denote the normal subgroup of $G(K)$ consisting of elements that act trivially on $V(\fZ_K)$. This subgroup has the following description:$$G(K)^{\natural} = \{g\in G(K)\,|\, \chi(g)\in \cO^{\times} \ {\rm {for\:all}} \  \chi\in {\rm{X}}^{\ast}_K(G_K)\}.$$  $G(K)^{\natural}$ contains $G'(K)$ and also every bounded subgroup of $G(K)$. Given a nonempty bounded subset $\Omega$ of an apartment $A$ of $\cB(G'/K)$, let $G(K)^{\Omega}$ denote the subgroup of $G(K)^{\natural}$ consisting of elements that fix $\Omega$ pointwise.  There is a smooth affine $\cO$-group scheme $\sG_{\Omega}$, with generic fiber $G_K$, whose group of $\cO$-rational points  considered as a subgroup of $G(K)$ is $G(K)^{\Omega}$ (when $G_K$ is quasi-split, these group schemes have been constructed in [BrT2,\:\S4]; for a simpler treatment of the existence and smoothness of these ``Bruhat-Tits group schemes'', see [Y]).  
 \vskip.5mm

 The subgroup $G(K)^{\Omega}$ is of finite index in the stabilizer of $\Omega$ in $G(K)^{\natural}$. In fact, any element of $G(K)^{\natural}$ which stabilizes $\Omega$ permutes the facets of the building that meet $\Omega$, and hence a subgroup of finite index of the stabilizer of $\Omega$ (in $G(K)^{\natural}$) keeps each facet that meets $\Omega$ stable and fixes every vertex of such a facet, hence it fixes pointwise every facet that meets $\Omega$. Thus a subgroup of finite index of the stabilizer of $\Omega$ in $G(K)^{\natural}$ fixes $\Omega$ pointwise, therefore  this subgroup is contained in $G(K)^{\Omega}$.  As $G(K)^{\Omega}\,(=\sG_{\Omega}(\cO))$ is a bounded subgroup of $G(K)$, the stabilizer of $\Omega$ in $G(K)^{\natural}$ is bounded. 
 
 \vskip.5mm
 The neutral component $\sG^{\circ}_{\Omega}$ of $\sG_{\Omega}$ is by definition the union of the connected generic fiber $G_K$ and the identity component of the special fiber of $\sG_{\Omega}$. The neutral component  $\sG^{\circ}_{\Omega}$  is an open $\cO$-subgroup scheme of $\sG_{\Omega}$, and it is affine [PY2, Lemma in \S3.5].  The subgroup $\sG^{\circ}_{\Omega}(\cO)$ is of finite index in $\sG_{\Omega}(\cO)$\,[EGA\,IV$_3$, Cor.\,15.6.5]. According to  [BrT2,\,1.7.1-1.7.2]  the $\cO$-group schemes $\sG_{\Omega}$ and $\sG^{\circ}_{\Omega}$ 
 are ``\'etoff\'e'' and hence by (ET) of [BrT2,\,1.7.1] their affine rings have the following description:$$\cO[\sG_{\Omega}] = \{ f\in K[G]\:\vert\:  f(\sG_{\Omega}(\cO))\subset \cO\};\,\, \, \cO[\sG^{\circ}_{\Omega}] = \{ f\in K[G]\:\vert\: f(\sG^{\circ}_{\Omega}(\cO))\subset \cO\}.$$  
 
\vskip.5mm
  If the above apartment $A$ corresponds to the maximal $K$-split torus $T$ of $G_K$, then there is a closed $\cO$-torus $\sT$ in $\sG^{\circ}_{\Omega}$ with generic fiber $T$. The special fiber $\osT$ of $\sT$ is a maximal $\kappa_s$-torus in the special fiber $\osG^{\circ}_{\Omega}$ of $\sG^{\circ}_{\Omega}$. Note that $\sT(\cO)$ is the maximal bounded subgroup of $T(K)$. For $\Omega$ as above, $\sG^{\circ}_{\Omega} = \sG^{\circ}_{\overline{\Omega}}$.  If  $\Omega$ is a nonempty subset of a facet $F$ of $\cB(G'/K)$, then $\sG^{\circ}_{\Omega}= \sG^{\circ}_{F}$, and  if moreover, $G$ is semi-simple, simply connected and quasi-split over $K$, then $\sG^{\circ}_{\Omega}(\cO)$  is the stabilizer of $\Omega$  in $G(K)$, so $\sG_{\Omega} = \sG^{\circ}_{\Omega}$, i.e., both the fibers of $\sG_{\Omega}$ are connected. 
  \vskip1mm
  
  As usual,  $G(K)^+$ will denote the normal subgroup of $G(K)$ generated by the $K$-rational elements of the unipotent radicals of parabolic $K$-subgroups of $G_K$. If $G$ is semi-simple and $\pi: \widehat{G}\rightarrow G$ is the simply connected central cover of $G$, then $\pi(\widehat{G}(K)^+) =G(K)^+$. So, if $G$ is semi-simple, simply connected and quasi-split over $K$\ and $\Omega$ is a subset of a facet of $\cB(G/K)$, the stabilizer of $\Omega$  in $G(K)^+$  fixes $\Omega$ pointwise. Moreover, assuming that $G$ is semi-simple and quasi-split over $K$, if an element  $g$ of $G(K)^+$ belongs to $\sG_{\Omega}(\cO)$, then it is actually contained in $\sG^{\circ}_{\Omega}(\cO)$.  
  \vskip1mm
  
  For a facet $F$ of $\cB(G'/K)$,  $\sG^{\circ}_F$ and $\sG^{\circ}_F(\cO)$ are respectively called  the {\it Bruhat-Tits parahoric group scheme} and the {\it parahoric subgroup} of $G(K)$ associated to  $F$. 
 The subset of points (of $\cB(G'/K)$) fixed under $\sG^{\circ}_F(\cO)$ is $\overline{F}$. 

 \vskip2mm

 \ni{\bf 1.9.}  We introduce the following partial order ``$\prec$'' on the set of nonempty subsets of $\cB(G'/K)$: Given two nonempty subsets  $\Omega$ and $\Omega'$ of $\cB(G'/K)$,  $\Omega'\prec \Omega$ if the closure $\overline\Omega$ of $\Omega$ contains $\Omega'$. For facets $F$ and $F'$ of $\cB(G'/K)$, if $F'\prec F$, we say that  $F'$ is a {\it face} of $F$. In a collection $\cC$ of facets, thus a facet is {\it maximal} if it is not a proper face of any facet belonging to  $\cC$, and a facet is {\it minimal} if no proper face of it belongs to $\cC$.  
\vskip.5mm

Given nonempty bounded subset $\Omega$  and $\Omega'$ of an apartment of $\cB(G'/K)$, with $\Omega'\prec \Omega$,  the inclusion $G(K)^{\Omega}\subset G(K)^{\Omega'}$ gives rise to a $\cO$-group scheme homomorphism $\sG_{\Omega}\rightarrow \sG_{\Omega'}$ that is the identity homomorphism on the generic fiber $G_K$. This homomorphism restricts to a $\cO$-group scheme homomorphism $\rho_{\Omega',\Omega}:\sG^{\circ}_{\Omega}\rightarrow \sG^{\circ}_{\Omega'}$ and induces a $\kappa_s$-homomorphism ${\overline{\rho}}_{\Omega',\Omega}: \osG^{\circ}_{\Omega}\rightarrow \osG^{\circ}_{\Omega'}$.  The restriction of\, ${\overline{\rho}}_{\Omega',\Omega}$ to any torus of $\osG^{\circ}_{\Omega}$ is an isomorphism onto a torus of $\osG^{\circ}_{\Omega'}$. In particular, if $F'\prec F$ are two facets of $\cB(G'/K)$, then there is a $\cO$-group scheme homomorphism $\rho_{F',F}: \sG^{\circ}_F\rightarrow \sG^{\circ}_{F'}$ that is the identity homomorphism on the generic fiber $G_K$.  We will assume in this paper that:

(1) The kernel of the induced homomorphism ${\overline{\rho}}_{F',F}: \osG^{\circ}_{F}\rightarrow \osG^{\circ}_{F'}$ is a smooth unipotent $\kappa_s$-subgroup of $\osG^{\circ}_F$.
\vskip.5mm

(2) The image $\fp(F'/F): = \overline{\rho}_{F',F}(\osG^{\circ}_F)$ is a pseudo-parabolic $\kappa_s$-subgroup of $\osG^{\circ}_{F'}$. 
\vskip.5mm

(3) Let $T$ be a maximal $K$-split torus of $G_K$ such that the apartment of $\cB(G'/K)$ corresponding to $T$ contains $F$. Let $\sT$ be the closed $\cO$-torus of $\sG^{\circ}_F$ with generic fiber $T$, and let $\osT$ be the special fiber of $\sT$.  We consider $\osT$ to be a maximal $\kappa_s$-torus of $\osG^{\circ}_{F}$, as well as of $\osG^{\circ}_{F'}$ (under the homomorphism $\overline{\rho}_{F',F}$), and also of the maximal pseudo-reductive quotient $\orG_{F'} := \osG^{\circ}_{F'}/\sR_{u,\kappa_s}(\osG^{\circ}_{F'})$ of $\osG^{\circ}_{F'}$.   Let $x$ be a point of $F'$ and $v$ be a vector in $V(T)$ such that $v+x$ is a point of $F$.  Then the nonzero weights of $\osT$ in the Lie algebra of the pseudo-parabolic $\kappa_s$-subgroup  $\fp(F'/F)/ \sR_{u,\kappa_s}(\osG^{\circ}_{F'})$ of $\orG_{F'}$ are the roots $a$ of $\orG_{F'}$ (with respect to $\osT$) such that $v(a)\geqslant 0$.  

\vskip.5mm

(4) The inverse image of the subgroup $\fp(F'/F)(\kappa_s)$ of $\osG^{\circ}_{F'}(\kappa_s)$, under the natural homomorphism $\pi_{F'}:\sG^{\circ}_{F'} (\cO)\rightarrow \osG^{\circ}_{F'}(\kappa_s)$  is the image $\rho_{F',F}(\sG^{\circ}_F(\cO))$  of $\sG^{\circ}_F(\cO)$ in $\sG^{\circ}_{F'}(\cO)$.  
\vskip.5mm

(5) $F\mapsto \fp(F'/F)$ is an order-preserving bijective map of the partially-ordered set $\{F\ |\ F'\prec F\}$ onto the set of pseudo-parabolic $\kappa_s$-subgroups of \,$\osG^{\circ}_{F'}$ partially-ordered by the opposite of inclusion. 
\vskip1mm

Note that (4) implies that the inverse image $P_{F'}^+$ under $\pi_{F'}$ of the normal subgroup $\sR_{u,\kappa_s}(\osG^{\circ}_{F'})(\kappa_s) \,(\subset \fp(F'/F)(\kappa_s))$ of \,$\osG^{\circ}_{F'}(\kappa_s)$  is contained in the image of $\sG^{\circ}_F(\cO)$ in $\sG^{\circ}_{F'}(\cO)$. So $P_{F'}^+$ fixes every facet $F$, $F'\prec F$, pointwise.  (5) implies that a facet $C$ of  $\cB(G'/K)$ is a chamber (i.e., it is a maximal facet) if and only if $\osG_C^{\circ}$ does not contain a proper pseudo-parabolic $\kappa_s$-subgroup, or, equivalently, the maximal pseudo-reductive quotient of $\osG_C^{\circ}$ is commutative\,[CGP, Lemma 2.2.3]. (When  $G_K$ is quasi-split, the above assertions are proved in [BrT2, Thm.\,4.6.33].) 
\vskip2mm

 \ni{\bf 1.10. Bruhat-Tits theory for the derived subgroup $G'$ and Levi subgroups.} 
 
 \vskip2mm
 
For a nonempty bounded subset $\Omega$ of an apartment of the building $\cB(G'/K)$ of $G'(K)$, let $\sG_{\Omega}$ be the smooth affine $\cO$-group scheme as in the preceding subsection. Then we call the neutral component of  the canonical smoothening [BLR,\,7.1,\,Thm.\,5] (see also [PY2, 3.2]) of the schematic closure of $G'$ in $\sG_{\Omega}$ the Bruhat-Tits $\cO$-group scheme associated to $\Omega$ and $G'$.  Its generic fiber is $G'$. It is easily seen that these $\cO$-group schemes  have the properties described in 1.8 and 1.9, and hence if  Bruhat-Tits theory is available for $G$ over $K$, then it is also available for the derived subgroup $G'$ over $K$.  

\vskip1mm
 
 Let $S\,(\subset G)$ be a $k$-split torus and $M:=Z(S)$ be the centralizer of $S$ in $G$. Then $M$ is a connected reductive $k$-subgroup and Bruhat-Tits theory is available for $M$, as well as for its derived subgroup $\sD(M) = (M,M)$, over $K$. In fact, the enlarged Bruhat-Tits building $\cB(M/K)$ of $M(K)$ is the union of apartments of $\cB(G/K)$ corresponding to maximal $K$-split tori of $G_K$ containing $S_K$.  \vskip.5mm
 
 Let $\cS$ be the maximal $k$-torus contained in the center of $M$ that splits over $K$, and let $V(\cS_K)= \bR\otimes_{\bZ}{\rm{Hom}}_K({\rm{GL}}_1, \cS_K)$. Then $V(\cS_K)$ operates on each apartment of  the enlarged building $\cB(M/K)$ by translation. The quotient of $\cB(M/K)$ by $V(\cS_K)$ is the Bruhat-Tits building of $M(K)$, as well as that of $\sD(M)(K)$.  Its apartments are the quotients of the apartments of $\cB(M/K)$ by $V(\cS_K)$.   For any nonempty bounded subset $\Omega$ of such an apartment, let $\sS$ be the schematic closure of $S_K$ in $\sG^{\circ}_{\Omega}$. Then $\sS$ is a $\cO$-torus in $\sG^{\circ}_{\Omega}$ with generic fiber $S_K$, and the Bruhat-Tits smooth affine $\cO$-group scheme, with generic fiber $M_K$, associated to $\Omega$ is  the centralizer $\sM^{\circ}_{\Omega}$  of $\sS$  in $\sG^{\circ}_{\Omega}$\,(smoothness of such centralizers is  known; see, for example, [${\rm{SGA3_{II}}}$, Exp.\,XI, Cor.\,5.3] or [CGP, Prop.\,A.8.10(2)]). The special fiber of $\sM^{\circ}_{\Omega}$ is obviously connected. 
 \vskip2mm

 \vskip2mm

\ni{\bf 1.11.}  As before, let $G(K)^{\natural}$ denote the normal subgroup of $G(K)$ consisting of elements that act trivially on $V(\fZ_K)$; $G(K)^{\natural}$ contains $G'(K)$ and also all bounded subgroups of $G(K)$. Let $T$ be a maximal $K$-split torus of $G_K$, $N(T)$ be its normalizer, and $Z(T)$ be its centralizer, in $G_K$. Let $A$ be the apartment of $\cB(G'/K)$ corresponding to $T$ and $C$ be a chamber in $A$. Then the stabilizer of $A$ in $G(K)$ is $\sN:=N(T)(K)$.  Let $\sI$ be the stabilizer of $C$ in $G(K)^{\natural}$; $\sI$ is a bounded subgroup of $G(K)$.  The maximal bounded subgroup $\sZ_b$ of $\sZ : = Z(T)(K)$ fixes $C$ pointwise (1.8) and hence it is contained in $\sI$. Let $\sN^{\natural} = \sN\cap G(K)^{\natural}$ and $\sZ^{\natural} = \sZ\cap G(K)^{\natural}$. As $G'(K)$ acts transitively on the set of ordered pairs consisting of an apartment  and a chamber  in it, and any two chambers of a building lie on an apartment, we conclude that $G(K)^{\natural} = \sI\sN^{\natural}\sI$.   The Weyl group $\sN/\sZ$ of $G_K$ is finite. We fix a finite subset $\cS$ of $\sN^{\natural}$ that maps onto $\sN^{\natural}/\sZ^{\natural}$. Then $G(K)^{\natural} = \sI \cS\sZ^{\natural}\sI$.   It is obvious that a subset $\cX$ of $G(K)^{\natural}$ is bounded if and only if $\sI \cX\sI$ is bounded, or, equivalently, if and only if there exists a bounded subset $\cY$ of  $\sZ^{\natural}$ such that $\cX\subset \sI\cS \cY\sI$.  
 \vskip1mm
 
 The subgroup $\sZ_b$ of $\sZ$ has the following description (Proposition 1.5): An element $z\in \sZ$ belongs to $\sZ_b$ if and only if for every  $K$-rational character $\chi$ of $Z(T)$, $\omega(\chi(z)) = 0$. We fix a basis $\{\chi_j\}_{j=1}^{\dim T}$  of the group of $K$-rational characters of $Z(T)$. Then  the map $z\mapsto (\omega(\chi_j(z)))$ provides an embedding of $\sZ/\sZ_b$ into $\bZ^{\dim T}$ and so a subset of $\sZ$ is bounded if and only if its image in $\sZ/\sZ_b$ is finite, or, equivalently, if and only if it is contained in the union of finitely many cosets of $\sZ_b$ in  $\sZ$.  Thus $\cX\,(\subset G(K)^{\natural})$ is bounded if and only if there exist a finite subset $\{n_i\}\subset \sN^{\natural}$ such that $\cX\subset \bigcup_i \sI n_i \sI$. 
 \vskip1mm
 
 Using these observations, we prove the following proposition.
 
 \vskip2mm
 
 \ni{\bf Proposition 1.12.} {\em A subset $\cX$ of $G(K)^{\natural}$ is bounded if and only if for every  $x\in \cB(G'/K)$ the set $\{g\cdot x\, \vert\, g\in \cX\}$ is of bounded diameter. }
 \vskip 1mm
 
 So  if a nonempty closed convex subset of $\cB(G'/K)$ is stable under the action of a bounded subgroup $\cG$ of $G(K)$, then by the Bruhat-Tits fixed 
 point theorem (Proposition 3.2.4 of [BrT1]) it contains a point fixed by $\cG$. 
 \vskip2mm
 
 \ni{\it Proof.} It is easy to see that to prove the proposition it suffices to prove that $\cX$ is bounded if and only if for some  $x\in \cB(G'/K)$, the set $\{g\cdot x\, \vert\, g\in \cX\}$ is of bounded diameter. We will now use the notation introduced in 1.11 and choose a  $x_0\in C$ fixed by $\sI$.  We have observed in 1.11 that $\cX$ is bounded if and only if there is a finite subset $\{n_i\}$ of $\sN^{\natural}$ such that $\cX\subset \bigcup_i \sI n_i \sI$.  Let $d$ be a $G(K)$-invariant metric on $\cB(G/K)$. Then for every $g\in \sI n_i \sI$, $d(g\cdot x_0, x_0) = d(n_i\cdot x_0, x_0)$. This implies that $\sI n_i\sI\cdot x_0$ is a subset of bounded diameter for each $i$, proving the proposition.\hfill$\Box$   
\vskip2mm

   \ni{\bf 1.13.}  We will assume throughout  that there is an action of $\Gamma = {\rm{Gal}}(K/k)$ on $\cB(G'/K)$ by polysimplicial isometries  such that  the orbit of every point under this action is finite,  and for all $g\in G(K)$, $x\in \cB(G'/K)$, $\gamma\in \Gamma$, we have $\gamma(g\cdot x) = \gamma(g)\cdot \gamma(x)$\, (cf.\,[BrT2, 4.2.12]).  Thus,  there is an action of $\Gamma\ltimes G(K)$ on $\cB(G/K)\,(= V(\fZ_K)\times \cB(G'/K))$. According to  the Bruhat-Tits fixed point theorem, $\cB(G/K)$ contains a point fixed under $\Gamma$. 
 
 \vskip.5mm
 For any apartment $A$ of $\cB(G/K)$ , and $\gamma\in \Gamma$, $\gamma(A)$ is an apartment and the action map $A\rightarrow \gamma(A)$ is affine.  Therefore, if $T$ is a $k$-torus of $G$ such that $T_K$ is a maximal $K$-split torus of $G_K$, then the apartment $A_T$ of $\cB(G/K)$ corresponding to $T_K$ is stable under the action of \,$\Gamma$,  and  $\Gamma$  acts on $A_T$  by affine transformations through a finite quotient. 
 \vskip.5mm
 
 Given a nonempty bounded subset $\Omega$ of an apartment of $\cB(G'/K)$ that is stable under the action of $\Gamma$, 
 $\sG_{\Omega}(\cO)$ is stable under $\Gamma$ and hence the affine ring $\cO[\sG_{\Omega}]\,(\subset K[G])$ of $\sG_{\Omega}$ is stable under the natural action of $\Gamma$ on $K[G]= K\otimes_k k[G]$. In such cases (i.e., when $\Omega$ is stable under the action of $\Gamma$), the $\cO$-group scheme $\sG_{\Omega}$, and so also its neutral component,  admit unique descents to a smooth affine $\fo$-group schemes with generic fiber $G$; the affine ring of these descents are $(\cO[\sG_{\Omega}])^{\Gamma}= \cO[\sG_{\Omega}]\cap k[G]$ and $(\cO[\sG^{\circ}_{\Omega}])^{\Gamma}= \cO[\sG^{\circ}_{\Omega}]\cap k[G]$; see [BLR, \S6.2, Ex.\,B]. As it is unlikely to cause confusion, in the sequel whenever $\Omega$ is stable under $\Gamma$,  we will use $\sG_{\Omega}$ and $\sG^{\circ}_{\Omega}$ to denote these smooth affine $\fo$-group schemes. We will denote the special fibers of $\fo$-group schemes $\sG_{\Omega}$ and $\sG^{\circ}_{\Omega}$ by $\osG_{\Omega}$ and $\osG^{\circ}_{\Omega}$ respectively. Both are smooth affine $\kappa$-groups and $\osG^{\circ}_{\Omega}$ is the identity component of $\osG_{\Omega}$.  The maximal pseudo-reductive quotient of   $\osG^{\circ}_{\Omega}$ will be denoted by $\orG_{\Omega}$.  
 
 \vskip.5mm
 
 For a point $x\in \cB(G/K)$ fixed under $\Gamma$, we will denote $\sG^{\circ}_{\{x\}}$, $\osG^{\circ}_{\{x\}}$ and $\orG_{\{x\}}$ by $\sG^{\circ}_x$, $\osG^{\circ}_x$ and $\orG_x$ respectively. By definition,  $\sG^{\circ}_x$ and $\sG^{\circ}_x(\fo)$ are respectively  the {\it Bruhat-Tits parahoric group scheme} and the {\it parahoric subgroup} of $G(k)$ associated to the point $x$. 

 \vskip1mm
 
 Let $T$ be a $k$-torus of $G$ such that $T_K$ is a maximal $K$-split torus of $G_K$. Let $\Omega$ be a nonempty bounded subset of the apartment of $\cB(G'/K)$ corresponding to $T_K$. We assume that $\Omega$ is stable under the action of $\Gamma$. Then the $\cO$-torus $\sT$  of 1.8 admits a unique descent to a closed $\fo$-torus of\, $\sG^{\circ}_{\Omega}$; in the sequel we will denote this $\fo$-torus also by $\sT$.  The generic fiber of $\sT$ is $T$, its special fiber $\osT$ is a maximal $\kappa$-torus of \,$\osG^{\circ}_{\Omega}$, and  $\sG^{\circ}_{\Omega}(\fo)\cap T(k)$ is the maximal bounded subgroup of $T(k)$.  If the $k$-torus $T$ contains a maximal $k$-split torus $S$ of $G$, then the generic fiber of the maximal $\fo$-split subtorus $\sS$ of  $\sT$ is $S$ and the special fiber $\osS\,(\subset \osT)$ of $\sS$ is a maximal $\kappa$-split torus of \,$\osG^{\circ}_{\Omega}$.

 \vskip2mm
 In view of the results on descent of $\cO$-group schemes described above, it is obvious that to establish descent of  Bruhat-Tits theory from $K$ to $k$ for $G$ it only needs to be shown that $\cB(G/K)^{\Gamma}= V(\fZ_K)^{\Gamma}\times \cB(G'/K)^{\Gamma}$ is an affine building, or, equivalently, $\cB(G'/K)^{\Gamma}$ is an affine building. We have observed in 1.10 that Bruhat-Tits theory is available for $G'$ over $K$. 
% So, in \S3, while proving that $\cB(G'/K)^{\Gamma}$ is an affine building we may (and will) assume that $G$ is a connected semi-simple $k$-group such that Bruhat-Tits theory is available for it over $K$.
\vskip2mm

 \ni{\bf 1.14.} Let $\cB = \cB(G'/K)^{\Gamma}$; $\cB$ is closed and convex  and is stable under the action of $G(k)$ on $\cB(G'/K)$.   We will show that $\cB$ carries a natural structure of a polysimplicial complex\,(3.2),  its facets (or polysimplices) being the intersections with $\cB$ of facets of $\cB(G'/K)$ that are stable under $\Gamma$,  and it is a ``thick'' affine building. The maximal facets (maximal in the ordering defined in 1.9) will be called {\it chambers} of $\cB$. We will prove that the dimension of $\cB$, and so that of any chamber of $\cB$, is $r := k$-rank $G$. The apartments of $\cB$ are, by definition, the polysimplicial subcomplexes which are intersections of {\it special $k$-apartments} of $\cB(G'/K)$ (see 1.15 below) with $\cB$. We will show that the apartments of $\cB$ are affine spaces of dimension $r$ and they are in bijective correspondence with maximal $k$-split tori of $G$. To show that $\cB$, considered as a polysimplicial complex, is a ``thick'' building, we will verify that the following four conditions defining a building in [T1, 3.1] (cf.\,also [R,\,Thm.\,3.11]) hold:
 \vskip1mm
 
 \ni(B1) $\cB$ is {\it thick}, that is, any  facet of codimension 1 (i.e., of dimension $r-1$) is a face of at least three chambers.
 \vskip.5mm
 
 \ni(B2) The apartments of $\cB$ are {\it thin chamber complexes}{\footnote{A polysimplicial complex $\Delta$ of dimension $r$ is called a {\it chamber complex} if  every facet of $\Delta$ is a face of a chamber (i.e., a facet of dimension $r$) and any two chambers of $\Delta$ can be joined by a {\it gallery} (see Proposition 3.5 for the definition). A chamber complex is {\it thin} if any polysimplex of codimension 1 is a face of exactly two chambers.}}.
 \vskip.5mm
 
 \ni(B3) Any two facets of $\cB$ lie on an apartment of $\cB$.
 
 \vskip.5mm
 
 \ni(B4) If facets $\cF_1$ and $\cF_2$ are contained in the intersection of two apartments $\cA$ and $\cA'$ of $\cB$, then there is a polysimplicial  isomorphism $\cA\rightarrow \cA'$ which fixes $\cF_1$ and $\cF_2$ pointwise.

\vskip2mm

\ni{\bf 1.15. Special $k$-tori and special $k$-apartments.} A {\it special $k$-torus} in $G$ is a  $k$-torus $T\,(\subset G)$ that contains a maximal $k$-split torus of $G$  and $T_K$ is a maximal $K$-split torus of $G_K$. 
The apartment in $\cB(G/K)$, or in $\cB(G'/K)$,  corresponding to $T_K$, for a special $k$-torus $T$,   will henceforth be called a  {\it special $k$-apartment} corresponding to the (special) $k$-torus $T$. According to \,[BrT2, Cor.\,5.1.12], if Bruhat-Tits theory is available for $G$ over $K$ (for example, if $G$ is quasi-split over $K$), then $G$ contains a special $k$-torus. As this is an important and very useful result, we will give its proof in the next section (see Proposition 2.3).
\vskip1mm

It is clear from the definition that every special $k$-apartment is stable under the action of the Galois group $\Gamma$. If $x\ne y$ are two points of a special $k$-apartment $A$ which are fixed under $\Gamma$, then the whole straight line in $A$ passing through $x$ and $y$ is pointwise fixed under $\Gamma$.  
\vskip1mm

Using the Bruhat-Tits fixed point theorem  we see that a facet of $\cB(G'/K)$  is $\Gamma$-stable if and only if it contains a point fixed under $\Gamma$, i.e., the facet meets $\cB=\cB(G'/K)^{\Gamma}$. A facet in the building $\cB(G'/K)$ that meets $\cB$ will be called a {\it $k$-facet}. 
\vskip2mm

\ni{\bf 1.16.} Let $X$ be a nonempty convex subset of $\cB(G'/K)$ and $\cC$ be the set of facets of $\cB(G'/K)$, or facets lying in a given apartment $A$ of this building, that meet $X$. Then it is easy to see (Proposition 9.2.5 (i),\,(ii),\:of [BrT1]) that all maximal facets in $\cC$ are of equal dimension. If $F$ is maximal among the facets lying in $A$ which meet $X$, then every facet contained in $A$ that meets $X$ is contained in the affine subspace $A_F$ of $A$ spanned by $F$. The dimension of $A_F$ is equal to $\dim (F)$; in particular, $F$ is an open subset of $A_F$.  Moreover, $A\cap X$ is contained in the affine subspace of $A$ spanned by $F\cap X$.  So for any facet $F'$ in $A$, $\dim (F\cap X) \geqslant \dim (F'\cap X)$. 
\vskip1mm

As $\cB$ is a nonempty convex subset of $\cB(G'/K)$, the above assertions hold for $X=\cB$. Maximal $k$-facets of $\cB(G'/K)$ will be called $k$-{\it chambers}.  The $k$-chambers are of equal dimension, and moreover, for any $k$-chamber $C$, $\cC := C^{\Gamma} =C\cap\cB$ is a chamber of $\cB$. Conversely, given a chamber $\cC$ of $\cB$, the unique facet $C$ of $\cB(G'/K)$ that contains $\cC$ is a $k$-chamber and $\cC = C\cap\cB$.   Note that a $k$-chamber may not be a chamber (i.e., it may not be a facet of $\cB(G'/K)$ of maximal dimension); see, however, Proposition 2.4.

 \vskip2mm

\ni{\bf 1.17.} Given a nonempty $\Gamma$-stable bounded subset $\Omega$ of an apartment of $\cB(G'/K)$ and a nonempty $\Gamma$-stable subset $\Omega'$ of $\overline\Omega$, the homomorphism $\rho_{\Omega',\Omega}$ described in 1.9 descends to a $\fo$-group scheme homomorphism  $\sG^{\circ}_{\Omega}\rightarrow \sG^{\circ}_{\Omega'}$ that is the identity homomorphism on the generic fiber $G$. We shall denote this $\fo$-homomorphism also by $\rho_{\Omega',\Omega}$; it induces a $\kappa$-homomorphism $\overline{\rho}_{\Omega',\Omega}: \osG^{\circ}_{\Omega}\rightarrow \osG^{\circ}_{\Omega'}$ between the special fibers. In particular, if $F'\prec F$ are two $k$-facets of $\cB(G'/K)$, then there is a $\fo$-group scheme homomorphism $\sG^{\circ}_{F}\rightarrow \sG^{\circ}_{F'}$ that is the identity homomorphism on the generic fiber $G$.  The image of the induced homomorphism $\osG^{\circ}_{F}\rightarrow \osG^{\circ}_{F'}$ is a pseudo-parabolic $\kappa$-subgroup $\fp(F'/F)$ of $\osG^{\circ}_{F'}$, and $F\mapsto \fp(F'/F)$ is an order-preserving bijective map of the partially-ordered set $\{F\ |\ F'\prec F\}$ onto the set of pseudo-parabolic $\kappa$-subgroups of $\osG^{\circ}_{F'}$ partially-ordered by opposite of inclusion\,(1.9). 
\vskip1mm

Thus, $F$ is a maximal $k$-facet (i.e., it is a $k$-chamber) if and only if $\fp(F'/F)$ is a minimal pseudo-parabolic $\kappa$-subgroups of $\osG^{\circ}_{F'}$. Now note that the projection map $\osG^{\circ}_{F'}\rightarrow \orG_{F'}$ induces an inclusion preserving bijective correspondence between the pseudo-parabolic $\kappa$-subgroups of $\osG^{\circ}_{F'}$ and the pseudo-parabolic $\kappa$-subgroups of its maximal pseudo-reductive quotient $\orG_{F'}$ [CGP, Prop.\,2.2.10].   Hence, a $k$-facet $C$ of $\cB(G/K)$ is a $k$-chamber if and only if the pseudo-reductive $\kappa$-group $\orG_C$ does not contain a proper pseudo-parabolic $\kappa$-subgroup, or, equivalently, this pseudo-reductive group contains a unique maximal $\kappa$-split torus (this torus is central so it is contained in every maximal torus of $\orG_C$) [CGP, Lemma 2.2.3].      

\vskip5mm

\ni{\bf 2. Nine basic propositions}

\vskip5mm
\ni{\bf Proposition 2.1.} {\em Let $\sG$ be a smooth affine $\fo$-group scheme and $\osG:=\sG_{\kappa}$ be its special fiber.
\vskip1mm

$(\rm{i})$ Let $\osT$ be a $\kappa$-torus in $\osG$. There exists a closed $\fo$-torus $\sT$ in $\sG$ whose special fiber is $\osT$.
\vskip.7mm

$(\rm{ii})$ Let $\sT$ and $\sT'$ be two closed  $\fo$-tori in $\sG$ such that there is an element $\overline{g}\in \sG(\kappa)$ that conjugates $\sT_{\kappa}$ onto $\sT'_{\kappa}$. There exists a $g\in \sG(\fo)$ lying over $\overline{g}$ that conjugates $\sT$ onto $\sT'$.
\vskip.7mm

$(\rm{iii})$ Let $\sT$ be a closed $\fo$-torus in $\sG$. Then the normalizer  $N_{\sG}(\sT)$ of $\sT$ in $\sG$ is a closed smooth $\fo$-subgroup scheme of $\sG$. In particular, the natural homomorphism $N_{\sG}(\sT)(\fo)\rightarrow N_{\sG}(\sT)(\kappa)$ is onto. } 
\vskip2mm

\ni{\it Remark.} The proof of assertion (i) of this proposition (and also that of the next proposition) is  essentially same as the proof of Proposition 5.1.10 of [BrT2]. In (i), since the special fiber of $\sT$ is $\osT$, the character groups of $\sT_{\cO}$ and $\osT_{\kappa_s}$  are isomorphic as $\Gamma$-modules, $\Gamma={\rm{Gal}}(K/k)= {\rm{Gal}}(\kappa_s/\kappa)$. In particular, $\sT$ is split if  $\osT$ is split.      
\vskip2mm

\ni{\it Proof of Proposition 2.1.}  (i) Let $X$ be the character group of $\osT_{\kappa_s}$ considered as a $\Gamma$-module under the natural action of $\Gamma$ and $\kappa_s[X]$ (resp.\,$\cO[X]$) be the group ring of $X$ with coefficients in $\kappa_s$ (resp.\,$\cO$). Then the affine ring of the $\kappa$-torus $\osT$ is $(\kappa_s[X])^{\Gamma}$. Let  $\sT$ be the $\fo$-torus whose affine ring is $(\cO[X])^{\Gamma}$. Then, clearly, the special fiber $\sT_{\kappa}$ of $\sT$ is isomorphic to $\osT$ and the  character group of $\sT_{\cO}$ is isomorphic as a $\Gamma$-module to $X$. We fix a $\kappa$-isomorphism ${\overline\iota}: \sT_{\kappa}\rightarrow \osT\,(\subset \osG)$ and view it as a closed immersion of $\sT_{\kappa}$ into $\osG$.   According to a result of Grothendieck\,[$\rm{SGA3_{II}}$, Exp.\,XI, 4.2], the homomorphism scheme ${\underline{\rm {Hom}}}_{{{\rm{Spec}}(\fo)}{{\textendash}\rm{gr}}}(\sT, \sG)$ is  representable by a smooth $\fo$-scheme $\mathscr{X}$. Clearly, ${\overline\iota}\in \mathscr{X}(\kappa)$.  Now since $\fo$ is Henselian, the natural map $\mathscr{X}(\fo)\rightarrow \mathscr{X}(\kappa)$ is surjective\,[EGA\,IV, 18.5.17], and hence there is a $\fo$-homomorphism $\iota:\sT\rightarrow  \sG$ lying over $\overline\iota$, i.e., $\iota_{\kappa} = \overline\iota$.  As $\overline{\iota}$ is a closed immersion, using [$\rm{SGA3_{II}}$, Exp.\,IX, 2.5 and 6.6] we see that $\iota$ is also a closed immersion.   We identify $\sT$ with a closed $\fo$-torus of $\sG$ in terms of $\iota$. Then the special fiber of $\sT$ is $\osT$. This proves assertion (i).  
\vskip1mm

(ii) The transporter scheme $\fT:={\rm{Transp}}_{\sG}(\sT, \sT')$,  consisting of points of the scheme $\sG$ that conjugate $\sT$ onto $\sT'$, is a closed smooth $\fo$-subscheme of $\sG$ (see [C, Prop.\,2.1.2] or  [SGA3$_{\rm{II}}$, Exp.\,XI, 2.4bis]). Let $\overline\fT$ be the special fiber of $\fT$. Then $\overline{g}$ belongs to ${\overline{\fT}}(\kappa)$. Now as $\fo$ is Henselian, the natural map $\fT(\fo)\rightarrow {\overline{\fT}}(\kappa)$ is surjective [EGA\,IV$_4$, 18.5.17]. Therefore, there exists a $g\in \fT(\fo)$ lying over $\overline{g}$.  This $g$ will conjugate $\sT$ onto $\sT'$.
\vskip1mm

(iii) In the proof of assertion (ii), by taking $\sT' = \sT$ we conlclude (iii).\hfill$\square$
\vskip2mm

\ni{\bf Proposition 2.2.} {\em Let $T$ be a maximal $K$-split torus of $G_K$ and $\Omega$ be a nonempty 
bounded subset of the apartment $A$ of $\cB(G'/K)$ corresponding to $T$. Let $T'$ be another maximal 
$K$-split torus of $G_K$ and $A'$ be the corresponding apartment of $\cB(G'/K)$. Then $\Omega$ is contained in 
$A'$ if and only if any of the following three equivalent conditions hold: 
\vskip1mm

{\rm{(i)}} There is an element $g\in \sG^{\circ}_{\Omega}(\cO)$ 
such that $T' = gTg^{-1}$. This element carries $A$ to $A'$ and fixes $\Omega$ pointwise. 
\vskip1mm

{\rm{(ii)}} $\sG^{\circ}_{\Omega}$ contains a closed $\cO$-torus with generic fiber $T'$.  
\vskip1mm

{\rm{(iii)}} $\sG^{\circ}_{\Omega}(\cO)\cap T'(K)$ is the maximal bounded subgroup of $T'(K)$.}

\vskip2mm

When $G_K$ is quasi-split, the first assertion of this proposition is [BrT2, Prop. 4.6.28(iii)]. The proof given below is different from the proof in [BrT2].  
\vskip2mm

\ni{\it Proof.} We will use the preceding proposition, with $\cO$ in place of $\fo$,  and denote 
$\sG^{\circ}_{\Omega}$ by $\sG$, and its special fiber by $\osG$, in this proof. Let $\sT$ 
be the closed $\cO$-torus of $\sG$ with generic fiber $T$.    If $\Omega$ is contained in $A'$, then 
$\sG$ contains a closed $\cO$-torus with generic fiber $T'$. Let us assume now that $\sG$ contains a closed $\cO$-torus $\sT'$ with generic fiber $T'$. As the residue field $\kappa_s$ of $\cO$ is separably 
closed, the special fibers $\osT$ and $\osT'$ of $\sT$ and $\sT'$ are maximal split tori of $\osG$, and hence there is an element $\overline{g}$ of $\osG(\kappa_s)$ that  conjugates $\osT$ onto $\osT'$ [CGP, Thm.\,C.2.3].  Now Proposition 2.1(ii) implies that there exists a $g\in\sG(\cO)$ lying over $\overline{g}$ that conjugates $\sT$ onto $\sT'$. This element fixes $\Omega$ pointwise and conjugates $T$ onto $T'$ and hence carries $A$ to $A'$.   Hence $\Omega$ is contained in $A'$. Conversely, if there is an element  $g\in \sG(\cO)$ 
such that $T' = gTg^{-1}$, then $\sT' := g\sT g^{-1}$ is a closed $\cO$-torus of $\sG$ with generic fiber $T'$,  
and $g$ carries $A$ to $A'$ fixing $\Omega$ pointwise. 

\vskip1mm

By Lemma 4.1 of [PY2], $\sG(\cO)\cap T'(K)$ is the maximal bounded subgroup of $T'(K)$ if and only if the schematic closure of $T'$ in $\sG$ is a $\cO$-torus. \hfill$\Box$

\vskip3mm

 \ni{\bf Proposition 2.3}\,([BrT2, Cor.\,5.1.12]){\bf .} {\em $G$ contains a special $k$-torus.}
\vskip2mm

\ni{\it Proof.} Let $S$ be a maximal $k$-split torus of $G$, $Z(S)$ its centralizer in $G$ and $Z(S)'$ be the derived subgroup of $Z(S)$. Then $Z(S)'$ is a connected semi-simple $k$-subgroup which is anisotropic over $k$. Let $S'$ be the maximal central $k$-torus of $Z(S)$ which splits over $K$. Then every special $k$-torus of $G$ that contains $S$ is of the form $S'\cdot T'$, where $T'$ is a special $k$-torus of $Z(S)'$. So after replacing $G$ with $Z(S)'$, we may (and we will) assume that $G$ is a semi-simple $k$-group that is anisotropic over $k$. 
\vskip1mm

Now let $x\in \cB(G'/K)$ be a point fixed under $\Gamma$\,(1.13). Let $\sG:= \sG_x^{\circ}$ be the smooth affine $\fo$-group scheme with generic fiber $G$ associated to $x$  in 1.13, and $\osG: =\sG_{\kappa}$ be the special fiber of $\sG$. Let $\osT$ be a maximal $\kappa$-torus of $\osG$.  
According to Proposition 2.1(i), there is a closed $\fo$-torus $\sT$ in $\sG$ with special fiber $\osT$.  Let $T$ be the generic fiber of $\sT$. Then $T$ is a $k$-torus of $G$ such that $T_K$ is a maximal $K$-split torus of $G_K$ since the special fiber $\osT$ of $\sT$ is a maximal $\kappa$-torus of $\osG$. Thus $T$ is a special $k$-torus of $G$. \hfill$\square$

%and $\sT$ be the $\fo$-torus whose character group over $\cO$ is isomorphic as a $\Gamma$-module to the character group of $\osT$ over $\kappa_s$. We fix a $\kappa$-isomorphism ${\overline\iota}: \sT_{\kappa}\rightarrow \osT$ and view it as an injective homomorphism of $\sT_{\kappa}$ into $\osG$.   According to a result of Grothendieck\,[$\rm{SGA3_{II}}$, Exp.\,XI, Cor.\,4.2], ${\underline{\rm {Hom}}}_{{{\rm{Spec}}(\fo)}{{\textendash}\rm{gr}}}(\sT, \sG)$ is representable by a smooth $\fo$-scheme $\mathscr{X}$. Clearly, ${\overline\iota}\in \mathscr{X}(\kappa)$.  Now since $\fo$ is Henselian, the natural map $\mathscr{X}(\fo)\rightarrow \mathscr{X}(\kappa)$ is surjective, and hence there is an injective $\fo$-homomorphism $\iota:\sT\rightarrow  \sG$ lying over $\overline\iota$.  The special fiber of the  image $\iota(\sT)$ is $\osT$ which is a maximal torus of $\osG$.  

\vskip2mm
\ni{\bf Proposition 2.4.} {\em  Every special $k$-apartment of $\cB(G'/K)$ contains a $k$-chamber. If $\kappa$ is perfect and of dimension $\leqslant 1$, then every $k$-chamber is a chamber of  $\cB(G'/K)$.}
\vskip2mm

\ni{\it Proof.} Let $A$ be a special $k$-apartment and $T$ be the corresponding special $k$-torus. Then $T$ contains a maximal $k$-split torus 
$S$ of $G$ and $T_K$ is a maximal $K$-split torus of $G_K$. As $A$ is stable under the action of $\Gamma$, by the Bruhat-Tits fixed point theorem, it contains a point $x$ which is fixed under $\Gamma$. Let $F$ be the facet lying on $A$ which contains $x$. Then, by definition, $F$ is a $k$-facet. Let $\sG^{\circ}_F$ be the  smooth affine $\fo$-group scheme, with connected fibers, associated to $\Omega=F$ in 1.13 and $\osG^{\circ}_F$ be the special fiber of $\sG^{\circ}_F$. Let $\sT$ be the closed $\fo$-torus of $\sG^{\circ}_F$ with generic fiber  $T$, and let $\sS$ be the maximal $\fo$-split subtorus of $\sT$ (cf.\,1.13). Then the generic fiber of $\sS$ is $S$. Let $\osS$ and $\osT$ be the special fibers of $\sS$ and $\sT$ respectively. We fix a minimal pseudo-parabolic $\kappa$-subgroup $\osP$ of $\osG^{\circ}_F$ containing $\osS$, then $\osP$ contains the centralizer of $\osS$ [CGP, Prop.\,C.2.4], and so it contains $\osT$.  Let $\cP$ be the inverse image of $\osP({\kappa_s})$ in $\sG^{\circ}_F(\cO)\, (\subset G(K))$ under the natural homomorphism $\sG^{\circ}_F(\cO)\rightarrow \osG^{\circ}_F({\kappa_s})$. Then $\cP$ is a parahoric subgroup of $G(K)$ contained in the parahoric subgroup  $\sG^{\circ}_F(\cO)$\: (see 1.17); $\cP$ contains $\sT(\cO)$ and is clearly stable under the action of $\Gamma$ on $G(K)$. Let $C$ be the facet of the Bruhat-Tits building $\cB(G'/K)$ fixed by $\cP$. Then $C$ contains $F$ in its closure and is stable under $\Gamma$, i.e., it is a $k$-facet; it is a $k$-chamber since $\osP$ is a minimal pseudo-parabolic $\kappa$-subgroup of $\osG^{\circ}_F$\,(1.17). Moreover, as $\cP$ contains the maximal bounded subgroup $\sT(\cO)$ of $T(K)$, $C$ lies in the apartment $A$\,(Proposition 2.2(iii)). 

If $\kappa$ is perfect and of dimension $\leqslant 1$, $\osG^{\circ}_F$ contains a Borel subgroup defined over $\kappa$ (1.7),  hence the minimal pseudo-parabolic subgroup $\osP$ is a Borel subgroup of $\osG^{\circ}_F$. So, in this case, $C$ is a chamber of the building $\cB(G/K)$. \hfill$\Box$    

\vskip2mm

\ni{\bf Remark 2.5.} Let $A$ be a special $k$-apartment of $\cB(G'/K)$. According to Proposition 2.4, there is a $k$-chamber contained in  $A$, so among the facets of $A$ that meet $\cB$, the maximal ones are $k$-chambers\,(1.16).  
\vskip2mm 

\ni{\bf Proposition 2.6.} {\em Given a $k$-chamber $C$ of the building $\cB(G'/K)$ that lies in a special $k$-apartment $A$, and a point $x\in \cB$, there is a special $k$-apartment that contains $C$ and $x$. Therefore, in particular, every point of $\cB$ lies in a special $k$-apartment.}
\vskip1.5mm

\ni{\it Proof.}    Let $T$ be the special $k$-torus corresponding to the apartment $A$.  Then $T$ contains a maximal $k$-split torus $S$ of $G$. We fix a point $y$ of $C\cap \cB$, then $\sG^{\circ}_y = \sG^{\circ}_C$. Let $\sS$ be the closed $\fo$-split torus in $\sG^{\circ}_C$ with generic fiber $S$. Let $\osS$ be the special fiber of $\sS$ and $\oS$ be the image of $\osS$ in $\orG_C$.  As $C$ is a $k$-chamber, $\oS$ is central and so every maximal torus of $\orG_y =\orG_C$  contains it (1.17). By the uniqueness of the geodesic $[xy]$, every point on it is fixed under $\Gamma$, i.e., $[xy] \subset \cB$. Restricted to any maximal $\kappa$-torus of $\osG^{\circ}_{[xy]}$, the composite $\kappa$-homomorphism  $$\overline{\rho}: \osG^{\circ}_{[xy]}\rightarrow \osG^{\circ}_y\rightarrow  \orG_y \,(=\orG_C),$$ where the first homomorphism is the $\kappa$-homomorphism $\overline{\rho}_{\Omega',\Omega}$ of 1.17  for $\Omega =[xy]$ and $\Omega' = \{y\}$,  is an isomorphism onto a maximal $\kappa$-torus of $\orG_C$. We fix a maximal $\kappa$-torus $\osT_{[xy]}$ of $\osG^{\circ}_{[xy]}$.   Then the maximal $\kappa$-split subtorus of $\osT_{[xy]}$ is isomorphic to $\oS$ since the isomorphic image ${\overline{\rho}}(\osT_{[xy]})$ of $\osT_{[xy]}$ is a maximal $\kappa$-torus of $\orG_C$ and so it contains $\oS$. 
\vskip1mm

According to Proposition 2.1(i), $\sG^{\circ}_{[xy]}$ contains a closed $\fo$-torus $\sT_{[xy]}$ whose special fiber (as a $\kappa$-subgroup of \,$\osG^{\circ}_{[xy]}$) is $\osT_{[xy]}$. The generic fiber $T_{[xy]}$ of $\sT_{[xy]}$ is then a $k$-torus of $G$ that splits over $K$ and contains a maximal $k$-split torus of $G$, so it is a special $k$-torus.  The special $k$-apartment of $\cB(G'/K)$ determined by $T_{[xy]}$ contains $[xy]$ and hence it contains $C$ and $x$.\hfill$\Box$ 
\vskip2mm
%
%\ni{\bf Proposition 2.6.} {\em Any $k$-facet is a face of a $k$-chamber.} 
%\vskip1.5mm
%
%\ni{\it Proof.} Let $F$ be a $k$-facet. Then, according to the preceding proposition, there is a special $k$-apartment $A$ containing it, and according to Proposition 2.4, $A$ contains a $k$-chamber.  So among the facets in $A$ that meet $\cB$ the maximal ones are $k$-chambers\,(1.16).   

\ni{\bf Proposition 2.7.} {\em Given points $x$, $y$ of $\cB$, there is a special $k$-apartment in $\cB(G'/K)$ that contains both $x$ and $y$. Therefore, given any two $k$-facets  of $\cB(G'/K)$ (which may not be different), there is a special $k$-apartment containing them.}

\vskip2mm

\ni{\it Proof.}  Let $F$ be the $k$-facet of $\cB(G'/K)$ that contains the point $y$. Let $C$ be a $k$-facet that contains $F$ in its closure, meets $\cB$, and is maximal among the facets with these two properties. Then $C$ is a $k$-chamber (Remark 2.5). Let $z\in C\cap \cB$. Then according to the previous proposition there is a special $k$-apartment which contains $z$, and hence also $C$. Now the same proposition implies that there is a  special $k$-apartment  which contains $C$ and $x$. This apartment then contains $\oC$, and hence also $F$, and so it contains $y$.\hfill$\Box$     
\vskip2mm

\ni{\bf Proposition 2.8.} {\em If $G$ is anisotropic over $k$, then $\cB= \cB(G'/K)^{\Gamma}\,(= \cB(G/K)^{\Gamma})$ consists of a single point.}
 \vskip2mm

 \ni{\it Proof.} To prove the proposition we will use Proposition 2.7. If $\cB$ contains points $x\ne y$, then according to that proposition there is a special $k$-apartment $A$ of  $\cB(G'/K)$  which contains both $x$ and $y$. Let $T$ be the special $k$-torus corresponding to $A$. Then $A$ is the image of  $V(T_K)+ x\,(\subset \cB(G/K))$ in $\cB(G'/K)$, where $V(T_K)= \bR\otimes_{\bZ}{\rm{X}}_{\ast}(T_K)$ with ${\rm{X}}_{\ast}(T_K) =  {\rm{Hom}}_K ({\rm{GL}}_1, T_K)$. As $G$ is anisotropic over $k$, $T$ is $k$-anisotropic, so ${\rm{X}}_{\ast}(T_K)^{\Gamma} = {\rm{Hom}}_k ({\rm{GL}}_1, T)$ is trivial. Hence, $A^{\Gamma}$ consists of a single point $x$. A contradiction!\hfill$\Box$ 
\vskip2mm

 Let $S$ be a maximal $k$-split torus of $G$. Let $Z(S)$ and $Z'(S)$ be the centralizers of $S$ in $G$ and $G'$ respectively, and $C$ be the central torus of $Z'(S)$.  Then $Z'(S)$ is an almost  direct product of $C$ and its derived subgroup $Z(S)'$. Let $\cC$ be the maximal $k$-subtorus of $C$ that splits over $K$, and $V(\cC_K) = \bR\otimes_{\bZ}{\rm{X}}_{\ast}(\cC_K)$. The {enlarged} Bruhat-Tits building $\cB(Z'(S)/K)= V(\cC_K)\times\cB(Z(S)'/K)$ of $Z'(S)(K)$ will be viewed as the union of apartments in the building $\cB(G'/K)$ corresponding to maximal $K$-split tori of $G_K$ that contain $S_K$. (Every such torus contains $\cC_K$. So $V(\cC_K)$ acts by translation on each apartment contained in  $\cB(Z'(S)/K)$.)    Let $S'\,(\subset \cC)$ be the maximal $k$-split subtorus of $C$. Then  $S'$ is the unique maximal $k$-split torus of $G'$ contained in $S$ and  $V(\cC_K)^{\Gamma} = V(S') := \bR\otimes_{\bZ}{\rm{X}}_{\ast}(S')$. Since $Z(S)'$ is $k$-anisotropic, according to the previous proposition $\cB(Z(S)'/K)^{\Gamma}$ consists of a single point. So $\cB(Z'(S)/K)^{\Gamma} = V(\cC_K)^{\Gamma}\times \cB(Z(S)'/K)^{\Gamma}$ is an affine space under $V(S')$.  Let $T$ be a special $k$-torus of $G$ containing $S$, i.e., a special $k$-torus  of $Z(S)$\:(such tori exist, see 1.15 or Proposition 2.3). The special $k$-apartment $A$ of $\cB(G'/K)$ corresponding to $T$ is a special $k$-apartment of $\cB(Z'(S)/K)$ (and special $k$-apartments of $\cB(Z'(S)/K)$ are in bijective correspondence with special $k$-tori of $G$ containing $S$).  The apartment $A$ is stable under $\Gamma$, so it contains a point,  say $x$,  fixed under $\Gamma$, and $A^{\Gamma}$ equals the set $V(S')+ x$ of  translates of $x$ under $V(S')$.    As $A^{\Gamma}$ is contained in $\cB(Z'(S)/K)^{\Gamma}$, and the latter is an affine space under $V(S')$, we conclude that $\cB(Z'(S)/K)^{\Gamma} = A^{\Gamma}$. We state this as the following proposition. 

\vskip2mm

\ni{\bf Proposition 2.9.} {\em For every special $k$-apartment $A$ of $\cB(Z'(S)/K)$, we have $A^{\Gamma} =\cB(Z'(S)/K)^{\Gamma}.$}
\vskip1mm

 $($Note that $A^{\Gamma}$  is an affine space under $V(S').)$}

\vskip2mm

\ni{\bf 2.10.} Let $N(S)$ be the normalizer of $S$ in $G$. As $N(S)$ normalizes $Z'(S)$, there is a natural action of  $N(S)(K)$ on $\cB(Z'(S)/K)$ and $N(S)(k)$ 
stabilizes $\cB(Z'(S)/K)^{\Gamma}$ under this action. For $n\in N(S)(K)$, the action of $n$ carries an apartment $A$ of $\cB(Z'(S)/K)$  to the apartment $n\cdot A$ by an affine transformation. 
\vskip1mm

 Now let $S'$ be the unique maximal $k$-split torus of $G'$ contained in $S$ and $T$ be a special $k$-torus of $G$ containing $S$.  Let $A:=A_T$ be the special $k$-apartment of $\cB(G'/K)$ corresponding to the torus $T$. As $T\supset S$, this apartment is contained in $\cB(Z'(S)/K)$ and it follows from the previous proposition that $\cB(Z'(S)/K)^{\Gamma} = A^{\Gamma}$. So we can view $\cB(Z'(S)/K)^{\Gamma}$ as an affine space under $V(S')=\bR\otimes_{\bZ}\X_{\ast}(S')$. We will now show, using the proof of the lemma in 1.6 of [PY1], that $\cB(Z'(S)/K)^{\Gamma}$ has the properties required of an apartment corresponding to the maximal $k$-split torus $S$ in the Bruhat-Tits building of $G(k)$ if such a building exists.  We need to check the following three conditions.  
 \vskip1.5mm
 
\ni{A1:} {\em The action of $N(S)(k)$ on $\cB(Z'(S)/K)^{\Gamma}= A^{\Gamma}$ is by affine transformations and the maximal bounded subgroup $Z(S)(k)_b$ of $Z(S)(k)$ acts trivially.}
\vskip1mm

Let ${\rm Aff}(A^{\Gamma})$ be the group of affine automorphisms of $A^{\Gamma}$ and $f: N(S)(k)\rightarrow {\rm Aff}(A^{\Gamma})$ be the action map.
\vskip1mm
 
\ni{A2:} {\em The group $Z(S)(k)$ acts on $\cB(Z'(S)/K)^{\Gamma}$ by translations, and the action is characterized by the following 
formula: for $z \in Z(S)(k)$,
$$\chi(f(z)) = -\omega(\chi(z)) \ {\rm for \ all}\  \chi \in {{\rm X}_k^{\ast}}(Z(S))\,(\hookrightarrow {{\rm X}_k^{\ast}}(S)),$$
here we regard the translation $f(z)$ as an element of  $V(S')\,(\hookrightarrow V(S)= \bR\otimes_{\bZ}{{\rm X}_{\ast}}(S))$}.
\vskip1mm

\ni{A3:} {\em For $g \in {\rm{Aff}}(A^{\Gamma})$, denote by $dg \in {\rm{GL}}(V(S'))$ the derivative of $g$. Then the map 
$N(S)(k)\rightarrow  {\rm{GL}}(V(S'))$, $n\mapsto df(n)$, is induced from the action of $N(S)(k)$ on ${{\rm X}_{\ast}}(S')$ $($i.e., it is the Weyl group action$)$.}
\vskip1mm

Moreover, as $G'$ is semi-simple, these three conditions determine the affine structure on $\cB(Z'(S)/K)^{\Gamma}$, see [T2, 1.2].
\vskip2mm

\ni{\bf Proposition 2.11.} {\em Conditions {\rm {A1, A2}} and {\rm{A3}} hold.}
\vskip2mm

\ni{\it Proof.} The action of $n \in N(S)(k)$ on $\cB(G'/K)$ carries the special $k$-apartment $A = A_T$  via an affine isomorphism $\varphi(n)\,:\, A\rightarrow A_{nTn^{-1}}$ to the special $k$-apartment $A_{nTn^{-1}}$ corresponding to the special $k$-torus $nTn^{-1}$ containing $S$. As $(A_{nTn^{-1}})^{\Gamma} = \cB(Z(S)/K)^{\Gamma} = A^{\Gamma}$, we see that $\varphi(n)$ keeps $A^{\Gamma}$ stable and so $\varphi(n)\vert_{A^{\Gamma}}$ is an affine automorphism of $A^{\Gamma}$.  
\vskip1mm

Let $V(T_K)= \bR\otimes_{\bZ}\,{\rm{Hom}}_K({\rm{GL}}_1, T_K)$ and  $V(nT_Kn^{-1})= \bR\otimes_{\bZ}\,{\rm{Hom}}_K({\rm{GL}}_1, nT_Kn^{-1})$. The derivative $d\varphi(n): V(T_K) \rightarrow V(nT_Kn^{-1})$  is induced from the map $${\rm {Hom}}_K({\rm {GL}}_1, T_K)={{\rm X}_{\ast}}(T_K)\rightarrow  {{\rm X}_{\ast}}(nT_Kn^{-1})={\rm{Hom}}_K({\rm{GL}}_1, nT_Kn^{-1}),$$ $\lambda\mapsto {\rm{Int}}\,n\cdot\lambda$, where ${\rm{Int}}\,n$ is the inner automorphism of $G$ determined by $n$. So, the restriction $df(n):V(S')\rightarrow V(S')$ is induced from the homomorphism ${\rm X}_{\ast}(S')\rightarrow {\rm X}_{\ast}(S')$, $\lambda\mapsto {\rm{Int}}\,n\cdot \lambda$. This proves {A3}.

\vskip1mm

Condition {A3} implies that $df$ is trivial on $Z(S)(k)$. Therefore, $Z(S)(k)$ acts on $\cB(Z'(S)/K)^{\Gamma} = A^{\Gamma}$ by translations. 
The action of the bounded subgroup $Z(S)(k)_b$ on $A^{\Gamma}$ admits a fixed point; see Proposition 1.12 and the observation 
following its statement. Hence, $Z(S)(k)_b$ acts by the trivial translation. This proves {A1}.
\vskip1mm

Since the image of $S(k)$ in $Z(S)(k)/Z(S)(k)_b\simeq \bZ^r$ is a subgroup of finite index, to prove the formula in {A2}, it suffices to prove it for $z\in S(k)$. But for $z\in S(k)$, $zTz^{-1} = T$, and $f(z)$ is a translation of the apartment $A$ ($f(z)$ is regarded as an element of $V(T_K)$) which satisfies\,(see 1.8):  $$\chi(f(z)) =-\omega(\chi(z))\ \ {\rm for \ all} \ \chi\in {\rm X}_K^{\ast}(T_K).$$ This implies the formula in {A2}, since the restriction map ${\rm X}_K^{\ast}(T_K)\rightarrow {\rm X}_K^{\ast}(S_K)\,(={\rm X}_k^{\ast}(S))$ is surjective and the image of the restriction map ${\rm X}_k^{\ast}(Z(S))\rightarrow {\rm X}_k^{\ast}(S)$ is of finite index in ${\rm X}_k^{\ast}(S)$. \hfill$\Box$
\vskip2mm

\ni{\bf 2.12.} By definition, the {\it apartments of $\cB$} are $A^{\Gamma}$, for special $k$-apartments $A$ of $\cB(G'/K)$. Let $T$ be a special $k$-torus of $G$, $S$ the maximal $k$-split torus of $G$ contained in $T$ and $S'$ be the maximal $k$-split torus of $G'$ contained in $S$. Let $A$ be the apartment of $\cB(G'/K)$ corresponding to $T$. Then (Proposition 2.9) $A^{\Gamma}=\cB(Z'(S)/K)^{\Gamma}$, thus the apartment $A^{\Gamma}$ of $\cB$  is uniquely determined by $S$,  and it is an affine space under $V(S')= \bR\otimes_{\bZ}{{\rm X}_{\ast}}(S')$. As maximal $k$-split tori of $G$ are conjugate to each other under $G(k)$, we conclude that $G(k)$ acts transitively on the set of apartments of $\cB$.   
\vskip5mm

\ni{\bf 3. Main results}
\vskip5mm

We will use the notations introduced in \S\S1,\,2. Thus $G$ will denote a connected reductive $k$-group, $G'$ its derived group. We assume that for $G$, and hence also for  $G'$, Bruhat-Tits theory is available over $K$.  As before, $\cB(G'/K)$ will denote the Bruhat-Tits building of $G(K)$; $\cB(G'/K)$ is also the Bruhat-Tits building of $G'(K)$. The enlarged Bruhat-Tits building $V(\fZ_K)\times \cB(G'/K)$ of $G(K)$ will be denoted by $\cB(G/K)$. In this section we will prove that $\cB(G'/K)^{\Gamma}$, and so also $\cB(G/K)^{\Gamma} = V(\fZ_K)^{\Gamma}\times \cB(G'/K)^{\Gamma}$, is a thick affine building. 
\vskip2mm
  
\ni{\bf Theorem 3.1.}  {\em Let $A_1$ and $A_2$ be special $k$-apartments of  $\cB(G'/K)$; $T_1$, $T_2$ be the corresponding special $k$-tori. Let $\Omega$ be a nonempty $\Gamma$-stable  bounded subset of $A_1\cap {A_2}$ and $\sG^{\circ}_{\Omega}$ be the smooth affine $\fo$-group scheme associated to $\Omega$ in 1.13. Then there is an element $g\in \sG^{\circ}_{\Omega}(\fo)\subset G(k)$ that carries ${A_1}^{\Gamma}$ onto ${A_2}^{\Gamma}$.
\vskip1mm

If the residue field $\kappa$  of $k$ is perfect and of dimension $\leqslant 1$, then there exists an element $g\in \sG^{\circ}_{\Omega}(\fo)\subset G(k)$ that  conjugates $T_1$ onto $T_2$, hence it carries the apartment $A_1$ onto the apartment $A_2$.
\vskip1mm

As $g$ belongs to $\sG^{\circ}_{\Omega}(\fo)$, it fixes $V(\fZ_K)\times\Omega$ pointwise.}

\vskip2mm

\ni{\it Proof.}  Let $S_1$ and $S_2$ be the maximal $k$-split tori of $G$ contained in $T_1$ and $T_2$ respectively. Let $\sG:= \sG^{\circ}_{\Omega}$ and $\sT_1$ and $\sT_2$ be the closed $\fo$-tori in $\sG$ with generic fibers $T_1$ and $T_2$ respectively (see 1.13). Let $\sS_1$ and $\sS_2$ be the maximal $\fo$-split subtori of $\sT_1$ and $\sT_2$ respectively.  Then the generic fibers of $\sS_1$ and $\sS_2$ are $S_1$ and $S_2$ respectively. Using Proposition 2.1(i) and the remark following that proposition, we see that the special fibers $\osS_1$ and $\osS_2$ of $\sS_1$ and $\sS_2$ respectively are maximal $\kappa$-split tori in the special fiber $\osG$ of $\sG$. Hence there exists an element $\overline{g}\in \osG(\kappa)$ that conjugates $\osS_1$ onto $\osS_2$ [CGP, Thm.\,C.2.3].  By Proposition 2.1(ii), there exists an element $g\in \sG(\fo)\,(\subset G(k))$  lying over $\overline{g}$ that conjugates $\sS_1$ onto $\sS_2$. As $g\sS_1 g^{-1} =\sS_2$, we infer that  $gS_1g^{-1} = S_2$, so $$g\cdot{A_1}^{\Gamma} =g\cdot\cB(Z'(S_1)/K)^{\Gamma} = \cB(Z'(S_2)/K)^{\Gamma} = {A_2}^{\Gamma},$$ and $g$ fixes $\Omega$ pointwise. 

\vskip1mm

To prove the second assertion of the theorem, let $\osT_1$ and $\osT_2$ be the special fibers of $\sT_1$ and $\sT_2$ respectively. Both of them are maximal $\kappa$-tori  of $\osG$. Now let us assume that $\kappa$ is perfect and of dimension $\leqslant 1$. Then the reductive  $\kappa$-group $\osG^{\rm{red}} := \osG/{\sR_{u,\kappa}}(\osG)$ is quasi-split\,(1.7)\,and hence any maximal $\kappa$-split torus of $\osG^{\rm{red}}$ is contained in a unique maximal torus.  Therefore, as the element $\overline{g}\in \osG(\kappa)$ chosen in the preceding paragraph conjugates $\osS_1$ onto $\osS_2$, it conjugates $\osT_1$ onto a maximal $\kappa$-torus of the solvable $\kappa$-subgroup $\osH: = \osT_2\cdot \sR_{u,\kappa}(\osG)$.  Since any two maximal $\kappa$-tori of the solvable $\kappa$-group $\osH$ are conjugate to each other under an element of $\osH(\kappa)$ [Bo, Thm.\,19.2], we conclude that $\osT_2$ is conjugate to $\osT_1$ under an element of $\osG(\kappa)$.  Now Proposition 2.1(ii) implies that there is an element $g\in \sG(\fo)\,(\subset (G(k))$ that conjugates $\sT_1$ onto $\sT_2$, so $gT_1g^{-1} = T_2$, and hence $g$ carries $A_1$ onto $A_2$ fixing $\Omega$ pointwise.\hfill$\Box$

\vskip3mm

\ni{\bf 3.2. Polysimplicial structure on $\cB=\cB(G'/K)^{\Gamma}$.} 
The {\it facets} (resp.\:{\it chambers}) of $\cB$ are by definition  the subsets $\cF:=F\cap\cB$ (resp.\:$\cC:=C\cap\cB$) for $k$-facets $F$ (resp.\:$k$-chambers $C$) of $\cB(G'/K)$. As the subset of points of $\cB(G'/K)$ fixed under $\sG^{\circ}_F(\cO) = \sG^{\circ}_{\cF}(\cO)$ is $\overline{F}$\,(1.8), the subset of points of $\cB$ fixed under $\sG^{\circ}_{\cF}(\cO)$ is $\overline{F}\cap \cB = \overline{\cF}$.  
\vskip1mm

Let $F$  be a minimal $k$-facet in $\cB(G'/K)$ and $A$ be a special $k$-apartment containing $F$\,(Proposition 2.6). We will presently show that  $F$ contains a unique point fixed under $\Gamma$ (i.e., $F$ meets $\cB$ in a single point).  Every special $k$-apartment is stable under the action of the Galois group $\Gamma$  which acts on it by affine automorphisms. Now if $x$ and $y$ are two distinct points in $F\cap \cB$, then the whole straight line in the apartment $A$ passing through $x$ and $y$  is pointwise fixed under $\Gamma$. This line must meet the boundary of $F$, contradicting the minimality of $F$. By definition, a {\it vertex} of $\cB$  is the unique point of $F\cap \cB$ for any minimal $k$-facet $F$ in $\cB(G'/K)$. 

\vskip1mm

Let $F$ be a $k$-facet in $\cB(G'/K)$ ($F$ is not assumed to be minimal) and $\sV_F$ be the set of vertices of $\cB$ contained in $\oF$. For $v\in\sV_F$, let $F_v$ be the face of $F$ which contains $v$. Since $v$ is a vertex of $\cB$, $F_v$ is a minimal $k$-facet. Now if $x$ and $y$ are two distinct vertices in $\sV_F$, then $\oF_x\cap \oF_y$ is empty. For this intersection is convex and stable under $\Gamma$ and hence if it is nonempty, it would contain a $\Gamma$-fixed point (i.e., a point of $\cB$). This would contradict the minimality of $k$-facets $F_x$ and $F_y$. Thus the sets of vertices (we call them $K$-vertices) of the facets $F_x$ and  $F_y$ are disjoint, and each one of these sets is $\Gamma$-stable. The union of the sets of $K$-vertices of $F_v$, for  $v \in\sV_F$, is the set of $K$-vertices of $F$. To see this, we observe that any $K$-vertex of $F$ is a $K$-vertex of a face of $F$ which is a minimal $k$-facet and so it contains  a (unique) point of $\sV_F$. Arguing by induction on dimension of $F$, we easily see that $\oF\cap \cB$ is the convex hull of the set $\sV_F$ of vertices of $\cB$ contained in $\oF$. The points of $\sV_F$ are by definition the vertices of the facet $\cF := F\cap \cB$ of $\cB$. 
\vskip1mm
 
 Given a $k$-facet $F$ of $\cB(G'/K)$, using the description of pseudo-parabolic $\kappa$-subgroups of  ${\orG_F}$ up to conjugacy, we see  (1.17)  that $\kappa$-rank of the derived subgroup of ${\orG_F}$ is equal to the codimension of $\cF :=F\cap \cB$ in $\cB$. 
 \vskip1mm
 
 Let $F$ be a $k$-facet of $\cB(G'/K)$, and $\cF = F\cap \cB$ be the corresponding facet of $\cB$. Then, for $g\in G(k)$, $g\cdot F$ is also a $k$-facet and $ g\cdot\cF= g\cdot(F\cap \cB) =(g\cdot F)\cap \cB$ is the facet of $\cB$ corresponding to $g\cdot F$. Thus the action of $G(k)$ on $\cB$ is by polysimplicial automorphisms. 
 
 \vskip1mm
  
We assume in this paragraph that $G$ is absolutely almost simple. Then the Bruhat-Tits building $\cB(G/K)$ is a simplicial complex, and in this case $\cB$ is also a simplicial complex with simplices $\cF:= F\cap \cB$, for $k$-facets $F$ of $\cB(G/K)$ ($F$ is a simplex!). To see this, note that given a nonempty subset ${\sV}'$ of $\sV_F$, the $k$-facet $F'$ whose set of $K$-vertices is the union of the set of $K$-vertices of $F_v$ for $v\in \sV'$ is a face of $F$, so $\cF': = F'\cap\cB$ is a face of $\cF$ and its set of vertices is $\sV'$.  
\vskip2mm

\ni{\bf 3.3.} If $G$ is semi-simple, simply connected and quasi-split over $K$, then for any $k$-facet $F$, the stabilizer of the facet $\cF = F\cap \cB$ of $\cB$ in $G(k)$ (resp.\,$G(K)$) is $\sG^{\circ}_{F}(\fo)$ (resp.\:$\sG^{\circ}_{F}(\cO)$), hence the stabilizer of $\cF$ fixes both $F$ and $\cF$ pointwise. This follows from the fact that  the stabilizer of $\cF$ also stabilizes $F$ since $F$ is the unique facet of $\cB(G/K)$ containing $\cF$. But, in case $G$ is semi-simple  simply connected and quasi-split over $K$, the stabilizer of $F$ in $G(K)$ is the subgroup $\sG^{\circ}_F(\cO)\,(=\sG^{\circ}_{\cF}(\cO)\subset G(K))$\,(1.8) and this subgroup fixes $F$ pointwise.     

\vskip2mm

\ni{\bf Proposition 3.4.} {\em Let $\cA$ be an apartment of $\cB$. Then there is a {unique} maximal $k$-split torus $S$ of $G$ such that $\cA = \cB(Z'(S)/K)^{\Gamma}$. So the stabilizer of $\cA$ in $G(k)$ is $N(S)(k)$.} 
\vskip2mm

\ni{\it Proof.}  We fix a maximal $k$-split torus of $G$ such that $\cA = \cB(Z'(S)/K)^{\Gamma}$. We will show that $S$ is uniquely determined by $\cA$. For this purpose, we observe that  as $N(S)(k)$ acts on $\cA$ and the maximal bounded subgroup $Z(S)(k)_b$ of $Z(S)(k)$ acts trivially (Proposition 2.11), the subgroup $\mathcal{Z}$ of $G(k)$ consisting of elements that fix $\cA$ pointwise is a bounded subgroup of $G(k)$ that is normalized by $N(S)(k)$ and contains $Z(S)(k)_b$. Using the Bruhat decomposition of $G(k)$ with respect to $S$, we see that every bounded subgroup of $G(k)$ that is normalized by $N(S)(k)$ is a normal subgroup of the latter. So the identity component of the Zariski-closure of $\mathcal{Z}$ is $Z(S)$. As $S$ is the unique maximal $k$-split torus of $G$ contained in $Z(S)$, both the assertions follow.\hfill$\Box$
\vskip2mm

\ni{\bf Proposition 3.5.} {\em Let $\cA$ be an apartment  of $\cB$, and $\cC$, $\cC'$ two chambers in $\cA$. Then there is a {gallery} joining $\cC$ and $\cC'$ in $\cA$, i.e., there is a finite sequence 
$$\cC=\cC_0, \:\cC_1, \:\ldots\,,\:\cC_m =\cC'$$ 
of chambers in $\cA$ such that for $i$ with $1\leqslant i\leqslant m$, $\cC_{i-1}$ and $\cC_i$ share a face of codimension 1.}
\vskip2mm

\ni{\it Proof.}  Let $\cA_2$ be the codimension $2$-skelton of $\cA$, i.e.,  the union of all facets in $\cA$ of codimension at least $2$. Then $\cA_2$ is a closed subset of $\cA$ of codimension $2$, so $\cA-\cA_2$ is a connected open subset of the affine space $\cA$. Hence $\cA-\cA_2$ is arcwise connected. This implies that given points $x\in \cC$ and $x'\in \cC'$, there is a piecewise linear curve in $\cA-\cA_2$ joining $x$ and $x'$. Now the chambers in $\cA$  that meet this curve make a gallery joining $\cC$ to $\cC'$. \hfill$\Box$
\vskip1mm

The dimension of any apartment, or any chamber, in $\cB$ is equal to the $k$-rank of $G' \,(= (G,G))$. 
A {\it panel} in $\cB$ is by definition a facet of codimension $1$.
\vskip2mm

\ni{\bf Proposition 3.6.} {\em Let $\cA$ be an apartment of  $\cB$ and $S$ be the maximal $k$-split torus of $G$ corresponding to this apartment. $($Then $\cA = (\cB(Z'(S)/K)^{\Gamma}.)$ The group $N(S)(k)$ acts transitively on the set of chambers of $\cA$.} 
\vskip2mm

\ni{\it Proof.} According to  the previous proposition, given any  two chambers in $\cA$, there exists a minimal gallery in $\cA$ joining these two chambers. So to prove the proposition by induction on the length of a minimal gallery joining two chambers, it suffices to prove that given two different chambers  $\cC$ and $\cC'$ in $\cA$ which share a panel $\cF$, there is an element $n\in N(S)(k)$ such that $n\cdot \cC = \cC'$. Let $\sG:=\sG^{\circ}_{\cF}$ be the smooth $\fo$-group scheme associated with the panel $\cF$ and $\sS\subset \sG$ be the closed $\fo$-split torus with generic fiber $S$. Let  $\osG$ be the special fiber of $\sG$, $\osS$ the special fiber of $\sS$. Then  $\osS$ is a maximal $\kappa$-split torus of $\osG$. The chambers $\cC$ and $\cC'$ correspond to minimal pseudo-parabolic $\kappa$-subgroups $\osP$ and $\osP'$ of $\osG$, see 1.17. Both of these minimal pseudo-parabolic $\kappa$-subgroups contain $\osS$ since the chambers $\cC$ and $\cC'$ lie on $\cA$. But then by Theorems C.2.5 and C.2.3 of [CGP], there is an element $\overline{n}\in \osG(\kappa)$ which normalizes $\osS$ and conjugates $\osP$ onto $\osP'$.   Now from Proposition 2.1(iii) we conclude that there is an element $n\in N_{\sG}(\sS)(\fo)$ lying over $\overline{n}$.  It is clear that $n$ normalizes $S$ and hence it lies in $N(S)(k)$; it fixes $\cF$ pointwise and $n\cdot\cC =\cC'$.\hfill$\Box$   
\vskip2mm

\ni{\bf Proposition 3.7.} {\em $\cB$ is thick, that is any panel is a face of at least three chambers, 
and every apartment of $\cB$ is thin, that is any panel lying in an apartment is a face of exactly two chambers of the apartment.}
\vskip2mm

\ni{\it Proof.} Let $F$ be a $k$-facet of $\cB(G'/K)$ that is not a chamber, and $C$ be a $k$-chamber of which $F$ 
is a face. Then there is an $\fo$-group scheme homomorphism $\sG^{\circ}_C\rightarrow \sG^{\circ}_F$. The image of $\osG^{\circ}_C$ in $\osG^{\circ}_F$,  under the induced homomorphism of special fibers, is a minimal pseudo-parabolic $\kappa$-subgroup of $\osG^{\circ}_F$, and conversely, any minimal pseudo-parabolic $\kappa$-subgroup of the latter determines a $k$-chamber with $F$ as a face. Now if $\kappa$ is infinite, $\osG^{\circ}_F$ clearly contains infinitely many minimal pseudo-parabolic $\kappa$-subgroups. On the other hand, if $\kappa$ is a finite field, then pseudo-parabolic $\kappa$-subgroups of  $\osG^{\circ}_F$ are parabolic and as any nontrivial irreducible projective $\kappa$-variety has at least three $\kappa$-rational points, we see that $F$ is a face of at least three  distinct $k$-chambers.
\vskip1mm

To prove the second assertion, let $\cF := F^{\Gamma}$ be a panel in an apartment $\cA$ of $\cB$, where $F$ is a $k$-facet in $\cB(G'/K)$. Let $S$ be the maximal $k$-split torus of $G$ corresponding to $\cA$. Let $\sG^{\circ}_F$ be the smooth affine $\fo$-group scheme associated with $F$ in 1.13 and $\orG_F$ be the maximal pseudo-reductive quotient of the special fiber of this group scheme. Let $\sS$ be the closed $\fo$-split torus of $\sG^{\circ}_F$ with generic fiber $S$. Then the chambers of $\cB$ lying in $\cA$ are in bijective correspondence with minimal pseudo-parabolic $\kappa$-subgroups of $\orG_F$ which contain the image $\oS$ of the special fiber of $\sS$\,(1.17). The $\kappa$-rank of the derived subgroup of $\orG_F$ is $1$ since $\cF$ is of codimension $1$ in $\cB$ (3.2). This implies that $\orG_F$ has exactly two minimal pseudo-parabolic $\kappa$-subgroups containing $\oS$. 

\vskip1mm
The second assertion also follows at once from the following well-known result in algebraic topology:\:In any simplicial complex whose geometric realization is a  topological manifold without boundary (such as an apartment $\cA$ in $\cB$), any simplex of codimension 1 is a face of exactly two chambers (i.e., maximal dimensional simplices).\hfill$\Box$  

\vskip2mm

We now assert that $\cB = \cB(G'/K)^{\Gamma}$ is an affine building. It is a polysimplicial complex (3.2).  Propositions 2.7, 3.5, 3.7 and Theorem 3.1 show that all the four conditions, recalled in 1.14, in the definition of buildings are satisfied for $\cB$, if $\cB(Z'(S)/K)^{\Gamma}\,(= \cB(Z'(S)/K)\cap \cB)$, for maximal $k$-split tori $S$ of $G$, are taken to be its apartments, and $\cF :=F\cap \cB$, for $k$-facets $F$ of $\cB(G'/K)$, are taken to be its facets. Thus we obtain the following:
\vskip2mm

\ni{\bf Theorem 3.8.} {\em $\cB=\cB(G'/K)^{\Gamma}$ is a building. Its apartments are the affine spaces 
$\cB(Z'(S)/K)^{\Gamma}$ under $V(S') :=\bR\otimes_{\bZ}{\rm{X}}_{\ast}(S')$, for maximal $k$-split tori $S'$ of $G' = (G,G)$.  Its chambers are $\cC := C\cap \cB$ for $k$-chambers $C$ of $\cB(G'/K)$, and its facets are $\cF:= F\cap \cB$ for $k$-facets $F$ of $\cB(G'/K)$. The group $G(k)$ acts on $\cB$ by polysimplicial isometries.}   

\vskip2mm

\ni {\bf Definition 3.9.}  $\cB$ is called the {\it Bruhat-Tits building of $G(k)$}.

\vskip2mm

 Since $G(k)$ acts transitively on the set of maximal $k$-split tori of $G$, it acts transitively on the set of apartments of $\cB$ (cf.\,2.12). Now  Proposition 3.6 implies the following:
 \vskip1.5mm
 
 \ni{\bf Proposition 3.10.} {\em G(k) acts transitively on the set of ordered pairs $(\cA, \cC)$ consisting of an apartment $\cA$ of $\cB$  and a chamber $\cC$ lying in the apartment $\cA$.} 
 \vskip2mm

\ni {\bf 3.11. Parahoric subgroups of $G(k)$.} For $x\in\cB$, the $\fo$-group scheme $\sG^{\circ}_x$ with connected fibers, described in 1.13, is by definition the Bruhat-Tits  parahoric $\fo$-group scheme,  and $\sG^{\circ}_x(\fo)$ is the parahoric subgroup of $G(k)$ associated to the point $x$. If $\cF$ is the facet of $\cB$ containing $x$ and $F$ is the $k$-facet of $\cB(G'/K)$ containing $\cF$, then $\sG^{\circ}_x = \sG^{\circ}_{\cF} = \sG^{\circ}_F$. The generic fiber of $\sG^{\circ}_x$ is $G$, and the subgroup $\sG^{\circ}_{x}(\fo)=\sG^{\circ}_F(\cO)^{\Gamma}\,(=\sG^{\circ}_F(\cO)\cap G(k))$ of $G(k)$ fixes $F$ pointwise. Since $F$ is the unique facet of $\cB(G'/K)$ containing $\cF$, the stabilizer of $\cF$ also stabilizes $F$. But $\sG^{\circ}_{F}(\cO)$ is of finite index in the stabilizer of $F$ in $G(K)$. Therefore, $\sG^{\circ}_x(\fo)=\sG^{\circ}_{\cF}(\fo)$ is of finite index in the stabilizer of $\cF$ in $G(k)$.   
For a $k$-chamber $C$ of $\cB(G'/K)$, let $\cC = C\cap \cB$ denote the  corresponding chamber of $\cB$. The  subgroup $\sG^{\circ}_{\cC}(\fo)$ is then a minimal parahoric subgroup of $G(k)$, and all minimal parahoric subgroups of $G(k)$ arise this way.  
\vskip.5mm

Let $P$ be a parahoric subgroup of $G(K)$ which is stable under the action of $\Gamma$ on $G(K)$, then the facet $F$ in $\cB(G'/K)$ corresponding to $P$ is $\Gamma$-stable, i.e., it is a $k$-facet. Let $\cF = F\cap \cB$ be the corresponding facet of $\cB$, and $\sG^{\circ}_{\cF}$ be the associated $\fo$-group scheme with generic fiber $G$ and with connected special fiber. Then $\sG^{\circ}_{\cF}(\fo) = \sG^{\circ}_F(\cO)^{\Gamma} = P^{\Gamma}$ is a parahoric subgroup of $G(k)$. Thus the parahoric subgroups of $G(k)$ are the subgroups of the form $P^{\Gamma}$, for $\Gamma$-stable parahoric subgroups $P$ of $G(K)$. 
\vskip2mm

\ni{\bf Proposition 3.12.} {\em The minimal parahoric subgroups of $G(k)$ are conjugate to each other under $G(k)$.}
\vskip 2mm

\ni{\it Proof.} The minimal parahoric subgroups of $G(k)$ are the subgroups $\sG^{\circ}_{\cC}(\fo)$ for chambers $\cC$ in the building $\cB$. Proposition 3.10  implies that $G(k)$ acts transitively on the set of chambers of $\cB$. \hfill$\Box$
\vskip2mm

\ni{\bf 3.13.} We say that $G$ is {\it residually quasi-split} if every $k$-chamber  in $\cB(G'/K)$  is actually a  chamber, or,  equivalently, if for any $k$-chamber $C$, the special fiber of the $\fo$-group scheme $\sG^{\circ}_C$ is solvable. If the residue field $\kappa$ of $k$ is perfect and of dimension $\leqslant 1$, then every semi-simple $k$-group is quasi-split over $K$ (1.7) and by Proposition 2.4, it is {residually quasi-split}. For residually quasi-split $G$, the minimal parahoric subgroups of $G(k)$ are called the {\it Iwahori subgroups} of $G(k)$. They are of the form $I^{\Gamma}$ for $\Gamma$-stable Iwahori subgroups $I$ of $G(K)$.       
\vskip2mm

\ni{\bf 3.14.} Assume that  $G$ is semi-simple, simply connected and quasi-split over $K$. Let  $\cF$ be a facet of $\cB$ and $F$ the $k$-facet of $\cB(G/K)$ containing $\cF$. 
Then the stabilizer of $\cF$ in $G(K)$ is $\sG^{\circ}_{\cF}(\cO)= \sG^{\circ}_F(\cO)$, so the stabilizer of $\cF$ in $G(k)$ is $\sG^{\circ}_{\cF}(\fo)$ and $\sG^{\circ}_{\cF} = \sG_{\cF} = \sG_F$, hence the stabilizer of $\cF$ in $G(k)$ fixes $\cF$ and $F$  pointwise\:(cf.\,1.8).  The normalizer of a parahoric subgroup $P$ of $G(k)$ is $P$ itself, for if $P$ is the stabilizer of the facet $\cF$ of $\cB$, then the normalizer of $P$ also stabilizes $\cF$, and hence it coincides with $P$.  

\vskip2mm

\ni{\bf 3.15. Tits systems in suitable subgroups of $G(k)$ provided by the building.}  We assume in this paragraph that $G$ is semi-simple.  Let $\cG$ be a subgroup of $G(k)$ that acts on $\cB$ by type-preserving automorphisms and acts transitively on the set of ordered pairs consisting of an apartment of $\cB$ and a chamber  lying in the apartment.  We fix an apartment $\cA$ of $\cB$ and a chamber $\cC$ lying in $\cA$. Let $S$ be the maximal $k$-split torus of $G$ corresponding to $\cA$ and $N(S)$ be the normalizer of $S$ in $G$. Let $B$ be the subgroup consisting of  elements in $\cG$  which stabilize  $\cC$, and $N$ be the group of elements in $\cG$ which stabilize $\cA$.   Then in view of  Theorem 3.8 and Proposition 3.10, according to [T1, Prop.\,3.11], $(B,N)$ is a saturated Tits system in $\cG$, and $\cB$ is the Tits building determined by this Tits system. Note that $\cG\cap \sG^{\circ}_{\cC}(\fo)$ is a subgroup of $B$ of finite index,  and $N = \cG\cap N(S)(k)$ since the stabilizer of  $\cA$ in $G(k)$ is $N(S)(k)$ by Proposition 3.4.

\vskip5mm  
\ni{\bf 4. Hyperspecial points of $\cB$ and hyperspecial parahoric subgroups of $G(k)$}
\vskip5mm

We will continue to use the notation introduced earlier.
\vskip1mm

\ni{\bf 4.1.} A point $x$ of $\cB$ is said to be a {\it hyperspecial point} if the $\fo$-group scheme $\sG^{\circ}_x$ is reductive. As the generic fiber $G$ of $\sG^{\circ}_x$ is reductive, the latter is reductive if and only if its special fiber $\osG^{\circ}_x$ is reductive. From the definition it is clear that every hyperspecial point of $\cB$ is also a hyperspecial point of the building $\cB(G/K)$ of $G(K)$. 
\vskip.5mm

In case $G$ is semi-simple, every hyperspecial point $x\in\cB$ is a vertex. In fact, if $\cF$ is the facet of $\cB$ containing $x$, and $y$ is a vertex of the compact polyhedron $\overline\cF$, then unless $x =y$, the image of the homomorphism $\overline{\rho}_{\{y\},\cF}: \osG^{\circ}_{\cF} = \osG^{\circ}_x\rightarrow \osG^{\circ}_y$, induced by the inclusion of $\{y\}$ in $\overline\cF$, is a proper pseudo-parabolic $\kappa$-subgroup of $\osG^{\circ}_y$.  However, as $\osG^{\circ}_{x}$ is reductive, its image  $\overline{\rho}_{\{y\},\cF}(\osG^{\circ}_{x})$ in $\osG^{\circ}_y$ is a reductive group. But a proper pseudo-parabolic subgroup cannot be reductive. We conclude that $x=y$, i.e., $x$ is a vertex. Moreover, since $x$ is a hyperspecial point of $\cB(G/K)$, it is also a vertex of this building.  
\vskip1mm

A {\it hyperspecial parahoric} subgroup of $G(k)$ is  by definition the parahoric subgroup  $\sG^{\circ}_x(\fo)$ for a hyperspecial point  $x$ of $\cB$.  Let $x\in \cB$ be a hyperspecial point and $A$ be a special $k$-apartment of $\cB(G/K)$ containing $x$ (Proposition 2.6). Let $T$ be the special $k$-torus corresponding to $A$, and $S$ be the maximal $k$-split torus of $G$ contained in $T$. Let $\sG^{\circ}_x$ be the reductive $\fo$-group scheme corresponding to $x$ and $\sS\subset \sT$ be the closed $\fo$-tori in $\sG^{\circ}_x$ with generic fibers $S\subset T$.  Let $Z_{\sG^{\circ}_x}(\sS)$ and $Z_{\sG^{\circ}_x}(\sT)$ respectively be the centralizers of $\sS$ and $\sT$ in $\sG^{\circ}_x$. Both these group subschemes are smooth (see, for example, [${\rm{SGA3_{II}}}$, Exp.\,XI, Cor.\,5.3] or [CGP, Prop.\,A.8.10(2)]), and hence their generic and special fibers are of equal dimension.  
\vskip.5mm

The special fibers $\osS$ and $\osT$ of $\sS$ and $\sT$ are respectively a maximal $\kappa$-split torus and a maximal $\kappa$-torus (containing $\osS$) of the special fiber $\osG^{\circ}_x$ of $\sG^{\circ}_x$\,(1.13).  As the residue field $\kappa_s$ of $K$ is separably closed, $\osT$ splits over $\kappa_s$ and hence the torus $T$  splits over $K$. Also, since $\osG^{\circ}_x$ is reductive, the centralizer of the maximal torus $\osT$ in $\osG^{\circ}_x$ is itself, so the special fiber of $Z_{\sG^{\circ}_x}(\sT)$ is $\osT$. By dimension consideration, this implies that $Z_{\sG^{\circ}_x}(\sT)= \sT$, so the centralizer of $T$ in $G$ equals $T$. Hence $T$ is a maximal torus of $G$ (and this maximal torus splits over $K$). Thus, if $\cB$ contains a hyperspecial point, then $G$ splits over $K$.      
\vskip.5mm

Now let us assume that the residue field $\kappa$ of $k$ is of dimension $\leqslant 1$. Then the reductive $\kappa$-group $\osG^{\circ}_x$ is quasi-split, i.e., it contains a Borel subgroup defined over $\kappa$\,(see 1.7), or, equivalently, the centralizer in $\osG^{\circ}_x$ of the maximal $\kappa$-split torus $\osS$ is a torus, and hence this centralizer is $\osT$. Thus the special  fiber of the group scheme  $(\sT\subset)\,Z_{\sG^{\circ}_x}(\sS)$ is $\osT$. From this we conclude that  $Z_{\sG^{\circ}_x}(\sS)= \sT$, so the centralizer of the maximal $k$-split torus $S$ of $G$ in the latter is the torus $T$. Therefore, $G$ is quasi-split.
\vskip.5mm

Thus we have proved the following:
\vskip1mm

\ni{\bf Proposition 4.2.}  {\em If the Bruhat-Tits building $\cB$ of $G(k)$ contains a hyperspecial point, then $G$ splits over the maximal unramified extension $K$ of $k$. Moreover, if the residue field $\kappa$ of $k$ is of dimension $\leqslant 1$, then $G$ is quasi-split over $k$.}   
\vskip2mm

Now we will establish the following partial converse of this proposition (cf.\,[BrT2, Prop.\,4.6.31]).
\vskip2mm

\ni{\bf Proposition 4.3.} {\em The Bruhat-Tits building of a quasi-split reductive $k$-group that splits over the maximal unramified extension $K$ of $k$ contains hyperspecial points.}  
\vskip2mm

\ni{\it Proof.} We begin by recalling a construction that produces all quasi-split reductive $k$-groups that split over $K$. Let $\sG$ be a Chevalley $\fo$-group scheme (i.e., a split reductive $\fo$-group scheme). Let $\sB$ be a Borel $\fo$-subgroup scheme of the adjoint group of $\sG$. Let $\sT$ be a maximal $\fo$-torus of $\sB$;  this torus splits over $\fo$. Let $G$, $B$ and $T$ be the generic fibers of  $\sG$, $\sB$ and $\sT$ respectively. Then $G$ is a $k$-split reductive group, $B$ is a Borel subgroup of the adjoint group of $G$ and $T$ is a $k$-split maximal torus of $B$. We may (and we will)  identify the outer automorphism group  ${\rm{Out}}(G/k)$ of $G/k$ with the subgroup of the automorphism group of $\sG$ that keeps $\sB$, $\sT$ and a pinning stable. Let $c$ be $1$-cocycle on ${\rm{Gal}}(K/k)$ with values in  ${\rm{Out}}(G/k)$. The Galois twist ${}_{c}G$ is a quasi-split $k$-group that splits over $K$, and all such groups arise in this way from suitable $\sG$ and $c$. Now we identify ${\rm{Gal}}(K/k)$ with the automorphism group of $\cO/\fo$ and consider $c$ to be a $1$-cocycle on this  automorphism group to obtain the twist ${}_c\sG$ of $\sG$. The generic fiber of the reductive $\fo$-group scheme ${}_c\sG$ is clearly ${}_c{G}$. 
\vskip.5mm

Thus, we conclude from the above that given a quasi-split reductive $k$-group that splits over $K$, there is reductive $\fo$-group scheme whose generic fiber is the given (quasi-split reductive) $k$-group. For simplicity, we now change notation. Let $G$ be a quasi-split reductive $k$-group that splits over $K$ and $\sG$ be a reductive $\fo$-group scheme with generic fiber $G$.  Let $\cB(G/K)$ and $\cB = \cB(G/K)^{\Gamma}$ be the Bruhat-Tits buildings of $G(K)$ and $G(k)$ respectively. Let $x$ be a point of $\cB(G/K)$ fixed by $\Gamma\ltimes \sG(\cO)$. The point $x$ lies in $\cB$. We will presently show that $x$ is a hyperspecial point of $\cB$.  

\vskip.5mm
As the smooth $\fo$-group scheme $\sG$ is ``\'etoff\'e'', the inclusion $\sG(\cO)\hookrightarrow \sG_x(\cO)$ induces a $\fo$-group scheme homomorphism $\sG\rightarrow \sG_x$ that is the identity homomorphism on the generic fiber $G$. Since the fibers of $\sG$ are connected, this homomorphism factors through a homomorphism $f: \sG\rightarrow \sG^{\circ}_x$. The kernel of the induced homomorphism $\overline{f}: \osG\rightarrow \osG^{\circ}_x$ between the special fibers is a unipotent normal subgroup scheme of the (connected) reductive group $\osG$, and so this kernel is zero-dimensional and  it follows by dimension consideration that $\overline{f}$ is surjective. Hence the special fiber $\osG^{\circ}_x$ of $\sG^{\circ}_x$ is reductive. This proves that $\sG^{\circ}_x$ is reductive. So, by definition, $x$ is a hyperspecial point of $\cB$.\hfill$\Box$   
\vskip2mm

\ni{\bf 4.4.} Let $x$ be a hyperspecial point of $\cB$. We will now show that the special fiber $\osG_x$ of $\sG_x$ is connected and hence $\sG_x = \sG^{\circ}_x$. We may (and do) replace $k$ with $K$ and assume that $\fo = \cO$ and the residue field $\kappa$ is separably closed.  Let $\sT$ be a maximal $\fo$-torus of $\sG^{\circ}_x$. Let $T$ be the generic fiber of $\sT$ and $\osT$ be its special fiber. Then $T$ is a maximal $k$-torus of $G$ and this torus is split\,(4.1). The centralizer $Z_{\sG_x}(\sT)$ of $\sT$ in $\sG_x$ is a smooth subgroup scheme ([${\rm{SGA3_{II}}}$, Exp.\,XI, Cor.\,5.3] or [CGP, Prop.\,A.8.10(2)]) containing $\sT$ as a closed subgroup scheme, and its generic fiber is $Z_G(T) = T$.  We see that the closed immersion $\sT\hookrightarrow Z_{\sG_x}(\sT)$ between smooth (and hence flat)  $\fo$-schemes is an equality on generic fibers and hence is an equality.  Therefore,  the inclusion $N_{\sG_x}(\sT)(\fo)\hookrightarrow N_G(T)(k)$ gives an embedding $N_{\sG_x}(\sT)(\fo)/\sT(\fo)\hookrightarrow N_G(T)(k)/T(k)$, and $Z_{\osG_x}(\osT) = \osT$. 
\vskip.5mm

As $\osT$ is a maximal $\kappa$-torus of $\osG_x$, by the conjugacy of maximal $\kappa$-tori in $\osG^{\circ}_x$ under $\osG^{\circ}_x(\kappa)$ ($\kappa$ is separably closed so every $\kappa$-torus is split!), we see that $\osG_x(\kappa) = N_{\osG_x}(\osT)(\kappa)\cdot \osG^{\circ}_x(\kappa)$. We know from Proposition 2.1(iii) that $N_{\sG_x}(\sT)(\fo)\rightarrow N_{\osG_x}(\osT)(\kappa)$ is surjective, and hence $N_{\sG_x}(\sT)(\fo)/\sT(\fo)\rightarrow N_{\osG_x}(\osT)(\kappa)/\osT(\kappa)$ is surjective too. So the order of $N_{\osG_x}(\osT)(\kappa)/\osT(\kappa)$ is less than or equal to that of $N_{\sG_x}(\sT)(\fo)/\sT(\fo)$ $(\hookrightarrow N_G(T)(k)/T(k))$.   On the other hand, $N_G(T)(k)/T(k)$ is the Weyl group of the root system of $(G,T)$ and $N_{\osG^{\circ}_x}(\osT)(\kappa)/\osT(\kappa)$ is the Weyl group of the root system of $(\osG^{\circ}_x,\osT)$, but these root systems are isomorphic\,([$\rm{SGA3_{III}}$, Exp.\,XXII, Prop.\,2.8]), hence their Weyl groups are isomorphic.    We conclude from these observations that the inclusion $N_{\osG^{\circ}_x}(\osT)(\kappa)/\osT(\kappa)\hookrightarrow N_{\osG_x}(\osT)(\kappa)/\osT(\kappa)$ is an isomorphism. So $N_{\osG_x}(\osT)(\kappa)=N_{\osG^{\circ}_x}(\osT)(\kappa)$, and therefore,  $$\osG_x(\kappa) = N_{\osG_x}(\osT)(\kappa)\cdot \osG^{\circ}_x(\kappa)= \osG^{\circ}_x(\kappa).$$ This implies that $\sG_x = \sG^{\circ}_x$.  
\vskip.5mm

\vskip5mm

\ni{\bf 5. Filtration of the root groups and valuation of root datum} 
\vskip5mm

\ni{\bf 5.1.} We fix a maximal $k$-split torus $S$ of $G$, and let $\Phi:=\Phi(G,S)$ be the root system of $G$ with respect to $S$. Let $\cB$ be the Bruhat-Tits building of $G(k)$ and $\cA$ be the apartment corresponding to $S$.  For a nondivisible root $a$, let $U_{a}$ be the root group corresponding to $a$. If $2a$ is also a root, the root group $U_{2a}$ is a subgroup of $U_a$. Let $S_a$ be the identity component of the kernel of $a$. Let $M_a$ be the centralizer of $S_a$ and $G_a$ be the derived subgroup of $M_a$. Then $M_a$ is a Levi-subgroup of $G$ and $G_a$ is a semi-simple subgroup of $k$-rank 1. Let $C_a$ be the central torus of $M_a$. Then $S_a$ is the maximal $k$-split subtorus of $C_a$. The root groups of $G_a$ and $M_a$ with respect to $S$ are $U_{\pm a}$, and also $U_{\pm 2a}$ in case $\pm2a$ are roots too. 
\vskip1mm

There is a $G_a(K)$-equivariant embedding of the  Bruhat-Tits building $\cB(G_a/K)$ of $G_a(K)$ into the Bruhat-Tits building $\cB(G/K)$ of $G(K)$, [BrT1, \S7.6];  such an embedding is unique up to translation by an element of $V((C_a)_K) := \bR\otimes_{\bZ} {\rm{X}}_{\ast}((C_a)_K)$.  Thus the set of $G_a(K)$-equivariant embeddings of $\cB(G_a/K)$ into $\cB(G/K)$ is an affine space under $V((C_a)_K)$ on which the Galois group $\Gamma$ of $K/k$ acts through a finite quotient. Therefore, there is a $\Gamma\ltimes G_a(K)$-equivariant embedding of  $\cB(G_a/K)$ into $\cB(G/K)$. This implies that there is a $G_a(k)$-equivariant embedding $\iota$ of the Bruhat-Tits building $\cB(G_a/K)^{\Gamma}$ of $G_a(k)$ into the Bruhat-Tits building $\cB(G/K)^{\Gamma}$ of $G(k)$. (In fact, such embeddings form an affine space under $V(S_a) := \bR\otimes_{\bZ}{\rm{X}}_{\ast}(S_a)= V((C_a)_K)^{\Gamma}$.) We shall consider the Bruhat-Tits building of $G_a(k)$, which is a Bruhat-Tits tree since $G_a$ is of $k$-rank 1, embedded in the Bruhat-Tits building of $G(k)$ in terms of $\iota$.  
\vskip2mm

\ni{\bf 5.2. Filtration of the root groups. } Given a real valued affine function $\psi$ on $\cA$ with gradient $a$, let $z$ be the point on the apartment $\cA_a$ ($\subset \cA$), corresponding to the maximal $k$-split torus of $G_a$ contained in $S$, in the Bruhat-Tits tree of $G_a(k)$, such that $\psi(z)=0$. Let $ \sG$ be the smooth affine $\fo$-group scheme with generic fiber $G_a$ such that $\sG(\cO)$ is the subgroup of $G_a(K)$ consisting of elements that fix $z$ (1.12 and 1.13). Let $\sG^{\circ}$ be the neutral component of $\sG$.  Let $\sS$ be the closed $1$-dimensional $\fo$-split torus of $\sG^{\circ}$ whose generic fiber is the maximal $k$-split torus of $G_a$ contained in $S$ and let  $\lambda: {\rm{GL}}_1\rightarrow \sS\,(\hookrightarrow \sG)$ be the $\fo$-isomorphism such that $\langle a,\lambda\rangle>0$. Let $\sU \,(:=U_{\sG}(\lambda)$, see [CGP, Lemma 2.1.5]) be the $\fo$-subgroup scheme of $\sG$ representing the functor $$R \leadsto \{ g\in \sG(R)\,|\, \lim_{t \to 0} \lambda(t)g\lambda(t)^{-1} = 1\}.$$ Using the last assertion of 2.1.8(3), and the first assertion of 2.1.8(4), of [CGP] (with $k$, which is an an arbitrary commutative ring in these assertions, replaced by $\fo$, and $G$ replaced by $\sG$), we see that $\sU$ is a closed smooth unipotent $\fo$-subgroup scheme with connected fibers (and hence $U_{\sG}(\lambda) = U_{\sG^{\circ}}(\lambda)$); its generic fiber is clearly $U_a$ and $\sU(\cO) = \sG(\cO)\cap U_a(K)$.  Since the fibers of $\sU$ are connected, it is actually contained in $\sG^{\circ}$. Denote by  $U_{\psi}$  the subgroup $\sU(\fo) =\sG(\fo)\cap U_a(k)\,(= \sG^{\circ}(\fo)\cap U_a(k))$ of $U_a(k)$.  For $\psi'\leqslant \psi$, $U_{\psi}\subseteq U_{\psi'}$ and the union of the $U_{\psi}$'s is $U_a(k)$. 
\vskip2mm

\ni{\bf 5.3.} We will now work with a given $u\in U_a(k)$--$\{1\}$.  Let $\psi_u$  be the largest real valued affine function on $\cA$ with gradient $a$ such that $u$ lies in $U_{\psi_u}$ and  let $z =z(u)$ be the unique point on the apartment $\cA_a$ where $\psi_u$ vanishes. We observe that $z$ is a  vertex in $\cA_a$. For otherwise, it would be  a point of a chamber $\cC$  (i.e., a 1-dimensional simplex) of $\cA_a$ and then since $u$ fixes $z$ it would fix the chamber $\cC$ pointwise, and hence it would fix both the vertices of $\cC$.  Now let $\psi> \psi_u$ be the affine function with gradient $a$ which vanishes at the vertex of $\cC$ where $\psi_u$ takes a negative value. Then $u$ belongs to $U_{\psi}$, contradicting the choice of $\psi_u$ to be the largest of such affine functions. As in the previous paragraph, let $\sG^{\circ}$ be the Bruhat-Tits group scheme, corresponding to the vertex $z = z(u)$,  with generic fiber $G_a$  and connected special fiber. Then $u$ lies in $\sG^{\circ}(\fo)$. Let $\osG^{\circ}$ be the special fiber of $\sG^{\circ}$. We assert that the image $\overline{u}$ of $u$ in $\osG^{\circ}(\kappa)$ does not lie in $\sR_{u,\kappa}(\osG^{\circ})(\kappa)$, for if it did, then $u$ would fix the unique chamber of $\cA_a$ which has $z$ as a vertex and  on which $\psi_u$ takes negative values. Then, as above, we would be able to  find an affine function $\psi>\psi_u$ with gradient $a$ such that $u\in U_{\psi}$, contradicting the choice of $\psi_u$.

\vskip1mm  

\ni {\bf 5.4.} Let $\sS$, and $\lambda: {\rm{GL}}_1\rightarrow \sS$ be as in 5.2 for $z = z(u)$. Let $\osS$ be the special fiber of $\sS$. Since $\osG^{\circ}(\kappa)$ contains an element  which normalizes $\osS$ and whose conjugation action on $\osS$ is by inversion, as in the proof of Proposition 3.6, by considering the smooth normalizer subgroup scheme $N_{\sG^{\circ}}(\sS)$, we conclude that $\sG^{\circ}(\fo)$ contains an element $n$ which normalizes $\sS$ and whose conjugation action on this torus is by inversion. 
\vskip1mm

 We shall now use the notation introduced in \S2.1 of [CGP].  According to Remark 2.1.11 and the last assertion of Proposition 2.1.8(3) of [CGP] (with $k$, which is an arbitrary commutative ring in that assertion, replaced by $\fo$, and $G$ replaced by $\sG^{\circ}$), the multiplication map $$U_{\sG^{\circ}}(-\lambda)\times Z_{\sG^{\circ}}(\lambda)\times U_{\sG^{\circ}}(\lambda)\rightarrow \sG^{\circ}$$ is an open immersion of $\fo$-schemes. We shall denote $U_{\sG^{\circ}}(\lambda)$, $Z_{\sG^{\circ}}(\lambda) (= Z_{\sG^{\circ}}(\sS))$ and $U_{\sG^{\circ}}(-\lambda)$ by $\sU_a$, $\sZ$ and $\sU_{-a}$ respectively, and the special fibers of these $\fo$-subgroup schemes by $\osU_a$, $\osZ$ and $\osU_{-a}$ respectively. Note that $\sU_{\pm a}$ are the $\pm a$-root groups of $\sG^{\circ}$ with respect to $\sS$, and $\osU_{\pm a}$ are the $\pm a$-root groups of $\osG^{\circ}$ with respect to $\osS$. Now since $n\sU_{-a}n^{-1} = \sU_{a}$, we see that $\Omega: = \sU_{-a}\sZ n\sU_{-a}$ is an open subscheme of $\sG^{\circ}$.  Let $\overline\Omega =\osU_{-a}\osZ n \osU_{-a} (\subset \osG^{\circ})$ be the special fiber of $\Omega$. 
\vskip1mm

Let $\pi: \osG^{\circ}\rightarrow \osG^{\circ}/\sR_{u,\kappa}(\osG^{\circ})$ be the maximal pseudo-reductive quotient of $\osG^{\circ}$. As $\overline{u}\notin \sR_{u,\kappa}(\osG^{\circ})(\kappa)$, $\pi(\overline{u})$ is a nontrivial element of $\pi(\osU_{a})(\kappa)$. Note that $\pi(\overline\Omega) = \pi(\osU_{-a})\pi(\osZ)\pi(n)\pi(\osU_{-a})$, and $\pi(\osU_{\pm a})$ are the $\pm a$-root groups of the pseudo-reductive $\kappa$-group $\osG^{\circ}/\sR_{u,\kappa}(\osG^{\circ})$ with respect to the maximal $\kappa$-split torus $\pi(\overline\sS)$ [CGP, Cor.\,2.1.9]. Now using Proposition C.2.24(i) of [CGP] we infer that $\pi(\overline{u})$ lies in $\pi(\overline\Omega)(\kappa)$.  We claim that  $\overline{u}\in \overline\Omega(\kappa)$. To establish this claim, it  would suffice to prove that $\overline{\Omega}\cdot{\sR_{u,\kappa}}(\osG^{\circ}) = \overline\Omega$.  
\vskip1mm

According to [CGP, Prop.\,2.1.12(1)],  the open immersion $$(\sR_{u,\kappa}(\osG^{\circ})\cap \osU_{a})\times (\sR_{u,\kappa}(\osG^{\circ})\cap\osZ)\times (\sR_{u,\kappa}(\osG^{\circ})\cap \osU_{-a})\rightarrow \sR_{u,\kappa}(\osG^{\circ}),$$ defined by multiplication, is an isomorphism of schemes. Using this, and the normality of $\sR_{u,\kappa}(\osG^{\circ})$ in $\osG^{\circ}$,  we see that $${\overline{\Omega}\cdot \sR_{u,\kappa}(\osG^{\circ})} = \osU_{-a}\osZ n \osU_{-a}\cdot \sR_{u,\kappa}(\osG^{\circ})= \osU_{-a}\osZ n\cdot\sR_{u,\kappa}(\osG^{\circ}) \osU_{-a}$$ $$= \osU_{-a}\osZ n(\sR_{u,\kappa}(\osG^{\circ})\cap \osU_{a})(\sR_{u,\kappa}(\osG^{\circ})\cap\osZ)\osU_{-a}=\osU_{-a}\osZ n \osU_{-a}= \overline\Omega.$$    

\vskip2mm
Now the following well-known lemma implies at once that $u$ is contained in $\Omega(\fo)$. Therefore, there exist $u', u''\in\sU_{-a}(\fo)$, such that $m(u):=u'uu''\in \sZ(\fo)n\,(\subset N_{\sG}(\sS)(\fo)\subset \sG^{\circ}(\fo))$. 
\vskip2mm
   
\ni{\bf Lemma 5.5.} {\em  Let $X$ be a scheme, and $\Omega \subset X$ an open subscheme.  If for a local ring $R$, $f: {\rm{Spec}}(R) \rightarrow X$ is a map carrying the closed point into $\Omega$, then $f$ factors through $\Omega$.}
\vskip2mm

\ni{\it Proof.}  Since $\Omega$ is an open subscheme of $X$, the property of $f$ factoring through $\Omega$ is purely topological; i.e., it is equivalent to show that the open subset $f^{-1}(\Omega) \subset {\rm{Spec}}(R)$ is the entire space.  Our hypothesis says that this latter open subset contains the closed point, so our task reduces to showing that the only open subset of a local scheme that contains the unique closed point is the entire space.  Said equivalently in terms of its closed complement, we want to show that the only closed subset $Z$ of ${\rm{Spec}}(R)$ not containing the closed point is the empty set. For an ideal $J \subset R$ defining $Z$, this is the obvious assertion that if $J$ is not contained in the unique maximal ideal of $R$ then $J = (1)$. \hfill$\Box$
\vskip2mm

\ni{\bf 5.6.} We recall that there exist unique $u', u''\in U_{-a}(k)$ such that $u'uu''$ normalizes $S$ [CGP, Prop.\,C.2.24(i)]. Thus the above $m(u)$ is uniquely determined by $u$. It acts on the apartment $\cA$ by an affine reflection $r(u)$ whose derivative (or, vector part) is the reflection associated with $a$. As $m(u)\in \sG^{\circ}(\fo)$, $r(u)$ fixes the point $z = z(u)$ defined above.  Hence, the fixed point set of the affine reflection  $r(u)$ is the hyperplane spanned by $S_a(k)\cdot{z}$ in $\cA$. As $\psi_u(z)=0$, this hyperplane is the vanishing hyperplane of the affine function $\psi_u$.   This observation implies at once that the filtration subgroups of $U_a(k)$ as defined in [T2, \S1.4] are same as the subgroups $U_{\psi}$ described above.  We also note that the largest half-apartment in $\cA$ that is fixed pointwise by the  element  $u$ is $\psi_u^{-1}([0, \infty))$. 
\vskip2mm

\ni{\bf 5.7.} As above, let $u', u''\in U_{-a}(k)$ be such that $m(u) = u'uu''$ normalizes $S$. Then $m(u) = (m(u)u''m(u)^{-1})u'u= uu''(m(u)^{-1}u'm(u))$. Since $m(u)^{-1}u'm(u)$ and $m(u)u''m(u)^{-1}$ belong to $U_a(k)$, we conclude that $m(u') = m(u) = m(u'')$. Hence, $\psi_{u'} = -\psi_u = \psi_{u''}$. Also, $m(u^{-1}) = m(u)^{-1}$, and hence $r(u) =r(u^{-1})$, and so $\psi_u = \psi_{u^{-1}}$.  
\vskip2mm

\ni{\bf 5.8.} Now assume that $2a$ is also a root of $G$ with respect to $S$, and $u\in U_{2a}(k)$--$\{1\}(\subset U_{a}(k)$--$\{1\})$. Let $u',\,u''$ be as in 5.6. Considering the semi-simple subgroup generated by the root groups $U_{\pm 2a}$, we see that $u',\,u''\in U_{-2a}(k)$. Let $\psi_u$ be the affine function as in 5.3. Then $2\psi_u$ is the affine function with gradient $2a$ whose vanishing hyperplane is the fixed point set of the reflection $r(u)$. Thus if we consider $u$ to be an element of $U_{2a}(k)$--$\{1\}$, then the associated affine function with gradient $2a$ is $2\psi_u$.  
\vskip2mm

\ni{\bf 5.9. Valuation of root datum.} The valuation $\varphi_a$ on the root group $U_a(k)$, corresponding to a given point $s\in \cA$ is defined as follows: For $u\in U_a(k)$--$\{1\}$, let $\varphi_a(u) = \psi_u(s)$.  According to a result of Tits (Theorem 10.11 of [R]), $(\varphi_a)_{a\in \Phi}$ is a valuation of the root groups $(U_a(k))_{a\in \Phi}$. From the results in 5.7, 5.8 we see that for $u\in U_a(k)$--$\{1\}$, if  $m(u) = u'uu''$ is as above, then $\varphi_{-a}(u')= -\varphi_a(u)= \varphi_{-a}(u'')$, and $\varphi_a(u) = \varphi_a(u^{-1})$. Moreover, if $2a$ is also a root, then $\varphi_{2a} = 2\varphi_a$ on $U_{2a}(k)$--$\{1\}$.
\vskip5mm

\ni{\bf 6. Residue field $\kappa$ perfect and of dimension $\leqslant 1$}

\vskip5mm

\ni{\em We will assume throughout this section that the residue field $\kappa$ of $k$ is perfect and is of dimension $\leqslant 1$. According to Proposition 2.4, then every $k$-chamber is a chamber of  $\cB(G/K)$; in other words,  every semi-simple $k$-group is residually quasi-split. } 
\vskip2mm

\ni{\bf Theorem 6.1.} {\em {\rm (i)} Any two special $k$-tori of $G$ are conjugate to each other under an element of $G(k)$.
\vskip1mm

 {\rm (ii)} Let $S$ be a maximal $k$-split torus of $G$, then any two special $k$-tori contained in $Z(S)$ are conjugate to each other under an element of the bounded subgroup $Z(S)'(k)$ of $Z(S)(k)$, where $Z(S)' =(Z(S), Z(S))$ is the derived subgroup of $Z(S)$.}
\vskip2mm

\ni{\it Proof.}  (i) For $i = 1, 2$, let $T_i$ be a special $k$-torus of $G$ and $A_i$ the corresponding special $k$-apartment of $\cB(G/K)$. If $A_1\cap A_2$ is nonempty, the first assertion follows immediately from the second assertion of Theorem 3.1. So let us assume that $A_1\cap A_2$ is empty. We fix a $k$-chamber $C_i$ in $A_i$, for $i=1, 2$ (Proposition 2.4). According to Proposition 2.7, there is a special $k$-apartment $A$ containing $C_1$ and $C_2$. Let $T$ be the special $k$-torus of $G$ corresponding to this apartment. Then using the second assertion of Theorem 3.1 twice, first for the pair $\{A, A_1\}$, and then for  the pair $\{A, A_2\}$ we see that $T$ is conjugate to both $T_1$ and $T_2$ under $G(k)$. So $T_1$ and $T_2$ are conjugate to each other under an element of $G(k)$.
\vskip1mm

(ii) Let $S'$ be the maximal central $k$-torus of $Z(S)$ which splits over $K$.  Then any special $k$-torus of $Z(S)$ is of the form $S'\cdot T'$, where $T'$ is a special $k$-torus of the semi-simple $k$-group $Z(S)'$. Now the second assertion follows from the first assertion applied to $Z(S)'$ in place of $G$. \hfill$\Box$ 

\hskip2mm

\ni{\bf Theorem 6.2.} {\em Let $T$ be a special $k$-torus of \,$G$ and $S$ be the maximal $k$-split torus of $G$ contained in $T$. Then $N(T)(k)\subset N(S)(k) = Z(S)'(k)\cdot N(T)(k)$. Therefore, the natural homomorphism $N(T)(k)\rightarrow N(S)(k)/Z(S)(k)_b$, induced by the inclusion of $N(T)(k)$ in $N(S)(k)$, is surjective.}
\vskip2mm

\ni{\it Proof.} Any $k$-automorphism of $T$ carries the unique maximal $k$-split subtorus $S$ to itself. So $N(T)(k)\subset N(S)(k)$.  Now let $n\in N(S)(k)$, then $nTn^{-1}$ is a special $k$-torus that contains $S$. So $T$ and $nTn^{-1}$ are special $k$-tori contained in $Z(S)$. Now Theorem 6.1(ii) implies that there is a $g\in Z(S)'(k)$ such that $g^{-1}Tg = nTn^{-1}$. Hence, $gn$ belongs to $N(T)(k)$, and $n =g^{-1}\cdot gn$. \hfill$\Box$   
\vskip3mm

The following result is in [BrT3, 4.4-4.5] for complete $k$.
\vskip2mm

\ni{\bf Theorem 6.3.} {\em Assume that $G$ is absolutely almost simple and anisotropic over $k$. Then it splits over the maximal unramified extension $K$ of $k$ and is of type $A_n$ for some $n$.}
\vskip2mm 

\ni{\it Proof.} We know from Proposition 2.8 that $\cB =\cB(G/K)^{\Gamma}$ consists of a single point, say $x$. Let $A$ be a special $k$-apartment of $\cB(G/K)$, and $C$ be a $k$-chamber in $A$ (Proposition 2.4). Then $C^{\Gamma} = C\cap \cB$ is nonempty, and hence it equals $\{x\}$.  Let $I$ be the Iwahori subgroup of $G(K)$ determined by the chamber $C$ and $T$ be the $k$-torus of $G$ corresponding to the apartment $A$. Then $I$ is stable under $\Gamma$, and $T_K$ is a maximal $K$-split torus of $G_K$. We consider the affine root system of  $G_K$ with respect to $T_K$ and let $\Delta$ denote its basis determined by the Iwahori subgroup $I$. Then $\Delta$ is stable under the natural action of $\Gamma$ on the affine root system and there is a natural $\Gamma$-equivariant bijective correspondence between the set of vertices of  $C$ and $\Delta$. As $\cB$ does not contain any facets of positive dimension, we see from the discussion in 3.2 that $\Gamma$ acts transitively on the set of vertices of $C$, and hence it acts transitively on $\Delta$. Now from the description of irreducible affine root systems, we see that $G_K$ is $K$-split and its root system with respect to the split maximal $K$-torus $T_K$ is of type $A_n$ for some $n$, for otherwise, the action of the automorphism group of the Dynkin diagram of $\Delta$ is not transitive on $\Delta$. \hfill$\Box$  
\vskip2mm

\ni{\bf Remark 6.4.}  If $k$ is a locally compact nonarchimedean field (that is, $k$ is complete and its residue field $\kappa$ is finite), then any absolutely almost simple $k$-anisotropic group $G$ is of {\it inner} type $A_n$ for some $n$. This assertion was 
proved by Martin Kneser for fields of characteristic zero, and Bruhat and Tits in general. In view of the previous theorem, to prove it, we just need to show that any simply connected absolutely almost simple $k$-group $G$ of {\it outer} type  $A_n$ for $n\geqslant 2$ is $k$-isotropic. Since there does not exist a noncommutative finite dimensional division algebra with center a quadratic Galois extension of $k$ which admits an involution of the second kind with fixed field $k$ (see [Sch, Ch.\,10, Thm.\,2.2(ii)]) if $G$ is of outer type, then there is a quadratic Galois extension $\ell$ of $k$ and a nondegenerate hermitian form $h$ on $\ell^{n+1}$ such that $G = {\rm{SU}}(h)$. But any hermitian form over a nonarchimedean locally compact field in at least 3 variables represents zero nontrivially, and hence ${\rm{SU}}(h)$ is isotropic for $n\geqslant 2$.   
\vskip1mm

The following example of an absolutely almost simple $k$-anisotropic group of {\it outer} type $A_{r-1}$ (over a discretely valued complete field $k$ with residue field of dimension $\leqslant 1$) was communicated 
to me by Philippe Gille. As usual, $\bC$ will denote the field of complex numbers; for a positive integer $r$, let $\mu_r$ denote the group of $r$-th roots of unity; $F = \bC(x)$ and $F{'}= \bC(x')$ with $x' = \sqrt{x}$. We take $k = F((t))$ and $k' = F'((t))$.  Since the Brauer groups of $F$ and $F'$ are trivial,  the residue maps  induce isomorphisms:
$${\rm{ker}}({}_r{\rm{Br}}(k' ) \xrightarrow[]{N_{k'/k}} {{}_r{\rm{Br}}}(k))\xrightarrow[]{\simeq} {\rm{ker}}({\rm{H}}^1 (F', \mu_r) \xrightarrow[]{N_{F'/F}} {\rm{H}}^1(F, \mu_r))$$ \vskip-1mm$\hskip6.7cm\xrightarrow[]{\simeq} {\rm{ker}}({F'}^{\times}/{{F'}^{\times}}^r \xrightarrow[]{N_{F'/F}} {F^{\times}}/{F^{\times}}^r);$

\ni see [S, \S2 of the Appendix after Ch.\,II]. The element $u := \frac{1+x'}{1-x'}\in {F'}^{\times}$   
has trivial norm over $F$, and has a pole of order $1$ at  $x' = 1$, so it cannot be an $r$-th power.
It defines a central simple $k'$-algebra $\sD$ which is division and cyclic of degree $r$. By Albert's theorem, $\sD$ carries a $k'/k$-involution $\tau$ of the second kind. The $k$-group ${\rm{SU}}(\sD,\tau)$ is of outer type $A_{r-1}$ and is anisotropic over $k$.  

\vskip2mm

In the following theorem, and in its proof,  we will use the notation introduced earlier in the paper, and assume, as before, that $k$ is a discretely valued field with Henselian valuation ring and perfect residue field $\kappa$ of dimension $\leqslant 1$.  
\vskip1mm

\ni{\bf Theorem 6.5.} {\em Let G be a simply connected semi-simple $k$-group. Then the Galois cohomology set ${\rm H}^1(k,G)$ is trivial.} 
\vskip2mm

\ni{\it Proof.} By Steinberg's theorem\,(1.7), ${\rm{H}}^1(K,G)$ is trivial, so ${\rm{H}}^1(k, G)\simeq {\rm{H}}^1(K/k, G(K))$.   Let $c: \gamma\mapsto c(\gamma)$ be a 1-cocycle on the Galois group $\Gamma$ of $K/k$ with values in $G(K)$ and $_cG$ be the Galois-twist of $G$ with the cocycle  $c$.  The $k$-groups $G$ and $_cG$ are isomorphic over $K$ and we will identify $_cG(K)$ with $G(K)$. (Recall that with identification of $_cG(K)$ with $G(K)$ as an abstract group, the ``twisted'' action of $\Gamma$ on $_cG(K)$ is described as follows: For $x\in {_cG}(K)$, and $\gamma\in \Gamma$, $\gamma\circ x= c(\gamma)\gamma(x)c(\gamma)^{-1}$, where $\gamma(x)$ denotes the $\gamma$-transform of $x$ considered as a $K$-rational element of the given $k$-group $G$.)   
\vskip1mm

Now let $I$ be a $\Gamma$-stable Iwahori subgroup of $G(K)$, say $I =\sG^{\circ}_C(\cO)$ for a $k$-chamber $C$ of $\cB(G/K)$. The subgroup $I$ is also an Iwahori subgroup of $_cG(K)$ (as $_cG(K)$ has been identified with $G(K)$ in terms of a $K$-isomorphism $({}_cG)_K\rightarrow G_K$). However, under the twisted action of $\Gamma$ on $_cG(K)$, $I$ may not be $\Gamma$-stable. But as $_cG$ is a residually quasi-split semi-simple $k$-group, $_cG(K)$ certainly contains an Iwahori subgroup which is stable under the twisted action of $\Gamma$. Since any two Iwahori subgroups of $_cG(K)$ are conjugate under $_cG(K) = G(K)$ (Proposition 3.10 for $K$ in place of $k$), there exists a $g\in G(K)$ such that $gIg^{-1}$ is stable under the twisted action of $\Gamma$. Then $c(\gamma)\gamma(g)I\gamma(g)^{-1}c({\gamma})^{-1} = gIg^{-1}$ for all $\gamma\in \Gamma$.   Hence, for $\gamma\in \Gamma$,  $c'(\gamma) := g^{-1}c(\gamma)\gamma(g)\in G(K)$ normalizes the Iwahori subgroup $I$.  As the normalizer of $I$ is $I$ itself (3.14 for $K$ in place of $k$), we conclude that $c'$, which is a $1$-cocycle on $\Gamma$ cohomologous to $c$, takes values in $I = \sG^{\circ}_C(\cO)$. So to prove the theorem, it suffices to prove the triviality of ${\rm{H}}^1(\Gamma, \sG^{\circ}_C(\cO))$. 
\vskip1mm

By unramified Galois descent over discrete valuation rings [BLR, \S6.2, Ex.\,B], this cohomology set classifies $\sG^{\circ}_C$-torsors $\mathscr{X}$ over $\fo$ which admit an $\cO$-point. (As $\mathscr{X}$ inherits $\fo$-smoothness from $\sG^{\circ}_C$, and $\cO$ is Henselian with algebraically closed residue field $\overline\kappa$, $\mathscr{X}$ does automatically admit $\cO$-points.)  Thus, it suffices to prove that every such torsor admits an $\fo$-point.  By $\fo$-smoothness of $\mathscr{X}$, and the henselian property of $\fo$, it suffices to prove that the special fiber of $\mathscr{X}$  admits a rational point. But the isomorphism class of the special fiber as a torsor is determined  by an element of the set ${\rm{H}}^1(\Gamma, \sG^{\circ}_C(\overline{\kappa}))$, and this cohomology set is trivial by Steinberg's Theorem\,(1.7) since $\kappa$ has been assumed to be perfect and of dimension $\leqslant 1$. \hfill$\Box$
%It follows from known results (see [SGA3$_{\rm{III}}$, Exp.\,XXIV, Prop.\,8.1(ii)] and [Gr, Thm.\,11.7 and Remark 11.8.3]) that the natural map ${\rm{H}}^1(\Gamma, \sG_C(\cO))\rightarrow {\rm{H}}^1(\Gamma, \sG_C({\overline{\kappa}}))$ is bijective. (We do not need surjectivity of this map which is more difficult to prove.) But since the special fiber of $\sG_C$ is a connected $\kappa$-group and $\kappa$ is of dimension $\leqslant 1$, by Steinberg's theorem (1.7) ${\rm{H}}^1(\Gamma, \sG_C({\overline{\kappa}}))$ is trivial. Hence, ${\rm{H}}^1(\Gamma, \sG_C(\cO))$ is also trivial. Therefore, $c'$ (hence also $c$) represents the trivial cohomology class. \hfill$\Box$   

 \vskip2mm

\ni{\bf Remark 6.6.} The above theorem was first proved by a case-by-case analysis by Martin Kneser for $k$ a nonarchimedean local field of characteristic zero with finite residue field. It was proved for all discretely valued {\it complete} fields $k$ with perfect residue field of dimension $\leqslant 1$ by Bruhat and Tits  [BrT3, Thm.\,in\:\S4.7].  If $\wk$ denotes the completion of $k$, then the natural map ${\rm{H}}^1(k,G)\rightarrow {\rm{H}}^1(\wk, G)$ is bijective [GGM, Prop.\,3.5.3(ii)]. So the vanishing theorem of Bruhat and Tits over the completion $\wk$ also implies the above theorem.

\vskip5mm

 \centerline{\bf References}
 
 \vskip4mm
 \ni[Be] V.\,G.\,Berkovich, {\it  \'Etale cohomology for non-archimedean analytic spaces.} 
 
 \ni Publ.\,Math.\,IHES {\bf 78} (1993), 5-161. 
 \vskip1.5mm

 \ni[Bo] A.\,Borel, {\it Linear algebraic groups} (second edition). Springer-Verlag, New York (1991).
 \vskip1.5mm
 
 \ni[BoT] A.\,Borel and J.\,Tits, {\it Groupes r\'eductifs}. Publ.\,Math.\,IHES {\bf 27} (1965), 55-150. 
 
 \vskip1.5mm
 
 \ni[BLR] S.\,Bosch, W.\,L\"utkebohmert and M.\,Raynaud, {\it N\'eron models}. Springer-Verlag, Heidelberg (1990).
\vskip1.5mm
 
 \ni[BrT1] F.\,Bruhat and J.\,Tits, {\it Groupes r\'eductifs sur un corps local.} Publ.\,Math.\,IHES {\bf 41}(1972).
 \vskip1.5mm
 
 \ni[BrT2] F.\,Bruhat and J.\,Tits, {\it Groupes r\'eductifs sur un corps local, II.} Publ.\,Math.\,IHES {\bf 60}(1984).
 \vskip1.5mm
 
  \ni[BrT3] F.\,Bruhat and J.\,Tits, {\it Groupes alg\'ebriques sur un corps local, III.} J.\,Fac.\,Sc.\, U.\,Tokyo {\bf 34}(1987), 671-698.
  \vskip1.5mm
   
   \ni[C]  B.\,Conrad, {\it Reductive group schemes}. In ``Group schemes:\:A celebration of SGA3'', Volume I, Panoramas et Synth\`eses \#{\bf 42-43}(2014), Soc.\,Math.\,France. 
\vskip1.5mm

\ni[CP] B.\,Conrad and G.\,Prasad, {\it Structure and classification of pseudo-reductive groups.} Proc.\,Symp.\,Pure Math.\,\#{\bf 94},  127-276, American Math.\,Soc.\,(2017). 

\vskip1.5mm

 \ni[CGP] B.\,Conrad, O.\,Gabber and G.\,Prasad, {\it Pseudo-reductive groups} (second edition). Cambridge U.\,Press, New York (2015).
 \vskip1.5mm

 \ni[$\rm{EGA\,IV_3}$] A.\,Grothendieck, {\it {\it El\'ements de g\'eom\'etrie alg\'ebrique}, IV: \'Etude locale des sch\'emas et des morphismes de sch\'emas},  Publ.\,Math.\,IHES {\bf 28}(1966),  5-255. 
 \vskip1.5mm 
 
 \ni[$\rm{EGA\,IV_4}$] A.\,Grothendieck, {\it {\it El\'ements de g\'eom\'etrie alg\'ebrique}, IV: \'Etude locale des sch\'emas et des morphismes de sch\'emas},  Publ.\,Math.\,IHES {\bf 32}(1967),  5-361. 
 \vskip1.5mm

 \ni[$\rm{SGA3_{II}}$] M.\,Demazure and A.\,Grothendieck, {\it Sch\'emas en groupes, Tome II}. Lecture Notes in Math.\,\#{\bf 152}, Springer-Verlag, Heidelberg (1970). 
 \vskip1.5mm
 
 \ni[$\rm{SGA3_{III}}$] M.\,Demazure and A.\,Grothendieck, {\it Sch\'emas en groupes, Tome III} (revised). Documents Math\'ematiques, Soci\'et\'e Math.\,France (2011). 
 \vskip1.5mm

\ni[GGM] O.\,Gabber, P.\,Gille and L.\,Moret-Bailly, {\it Fibr\'es principaux sur les corps valu\'es hens\'eliens.} Algebr.\,Geom.\,{\bf 1}(2014), 573-612. 
\vskip1.5mm
% 
% \ni[EGA] A.\,Grothendieck, {\it El\'ements de G\'eom\'etrie Alg\'ebrique}.
%Publ.\,Math.\,IHES {\bf 4, 8, 11, 17, 20, 24, 28, 32}, 1960--7.  
%\vskip1.5mm
%
%\ni[Gr] A.\,Grothendieck, {\it Groupe de Brauer III} in {\it Dix expos\'es sur la cohomologie des sch\'emas}, North-Holland, Amsterdam (1968). 
%\vskip1.5mm
  
\ni[P] G.\,Prasad, {\it Elementary proof of a theorem of Bruhat-Tits-Rousseau and of a theorem of Tits}. Bull.\,Soc.\,Math.\,France\,{\bf 110}(1982), 197-202.

\vskip1.5mm
 
 \ni[PY1] G.\,Prasad and J.-K.\,Yu, {\it On finite group actions on reductive groups and buildings.} Invent.\,Math.\,{\bf 147}(2002), 545--560.
\vskip1.5mm

\ni[PY2] G.\,Prasad and J.-K.\,Yu, {\it On quasi-reductive group schemes.} J.\,Alg.\,Geom. 
{\bf 15}\,(2006), 507-549.
\vskip1.5mm 

 \ni[R] M.\,Ronan, {\it Lectures on buildings}. University of Chicago Press, Chicago (2009).
 \vskip1.5mm

\ni[Sch] W.\,Scharlau, {\it Quadratic and hermitian forms}. Springer-Verlag, Heidelberg (1985).
\vskip1.5mm
  
\ni[S] J-P. Serre,  {\it Galois cohomology.} Springer-Verlag, New York (1997).
\vskip1.5mm

\ni[T1] J.\,Tits, {\it Buildings of spheirical type and finite {\rm BN}-pairs}.\:Lecture Notes in Math. \#{\bf 386}, Springer-Verlag, Berlin (1974).
\vskip1.5mm

\ni[T2] J.\,Tits, {\it Reductive groups over local fields.} Proc.\,Symp.\,Pure 
Math.\,\#{\bf 33}, Part I, 29--69, American Math.\,Soc.\,(1979).
\vskip1.5mm

\ni[Y] J.-K.\,Yu, {\it Smooth models associated to concave functions in Bruhat-Tits theory}. In ``Group schemes:\:A celebration of SGA3'', Volume III, Panoramas et Synth\`eses \#{\bf 47}(2016), Soc.\,Math.\,France. 

\vskip5mm

\ni{University of Michigan}

\ni{Ann Arbor, MI 48109.}

\ni{e-mail: gprasad@umich.edu}

\end{document}